\documentclass[10pt,twoside]{article}

\usepackage[centertags]{amsmath}
\usepackage{amsfonts}
\usepackage{amssymb}
\usepackage{amsthm}
\usepackage{epsfig}
\usepackage{color}
\allowdisplaybreaks[4]

\newcommand{\COM}[1]{}

\begin{document}
\baselineskip 15pt \setcounter{page}{1}
\title{\bf \Large  On the maxima of continuous and discrete time Gaussian order statistics processes \thanks{Research supported by National Science Foundation of China
(No. 11501250) and Natural Science Foundation of Zhejiang Province of China (No. LY18A010020).} }
\author{
{\small Zhongquan Tan}\footnote{ E-mail address:  tzq728@163.com }\\
{\small\it  College of Mathematics, Physics and Information Engineering, Jiaxing University, Jiaxing 314001, PR China;}\\
}
 \maketitle
 \baselineskip 15pt

\begin{quote}
{\bf Abstract:}\ \ In this paper, we study the asymptotic relation between the maximum of a
continuous order statistics process formed by stationary Gaussian processes
and the maximum of this process sampled at discrete time points. It is shown that,  these two maxima are
asymptotically independent when the Gaussian processes are weakly dependent and the discrete points are sufficient sparse,
while for other case,
 these two maxima are asymptotically dependent.

{\bf Key Words:}\ \ continuous time process,  discrete time process, extreme values,
order statistics processes.

{\bf AMS Classification:}\ \ Primary 60F05; secondary 60G15

\end{quote}

\section{Introduction}

Let $\{X(t), t\geq0\}$ be a stationary Gaussian process with mean 0,
variance 1, correlation function $r(t)$ and continuous sample
functions. Suppose the following conditions on
the correlation function $r(t)$ hold, i.e.,
for some $\alpha\in(0, 2]$,
\begin{eqnarray}
\label{eq1.1}
r(t)=1-|t|^{\alpha}+o(|t|^{\alpha})\ \ \mbox{as}\ \ t\rightarrow 0\ \ \mbox{and}\ \ r(t)<1\ \ \mbox{for}\ \ t>0
\end{eqnarray}
and
\begin{eqnarray}
\label{eq1.2}
r(t)\ln t\rightarrow r\in [0,\infty],\ \ \mbox{as}\ \  t\rightarrow \infty.
\end{eqnarray}
In the literature, the Gaussian process $\{X(t), t\geq0\}$ is called weakly dependent and strongly dependent if (\ref{eq1.2}) holds with $r=0$ and
$r\in (0,+\infty]$, respectively.

Let $(X_{1}(t),\ldots,X _{n} (t)), t\geq 0$ be a Gaussian vector process, the components of which are
independent copies of  the Gaussian process $\{X(t), t\ge 0\}$.
Let $\{X_{m:n}(t), t\geq 0\}$  with $1\leq m\leq n$, be the $m$th upper order statistic process  of $\{X_{i}(t), t\geq 0\}$, $1\leq i\leq n$,  defined by
\begin{eqnarray}
\label{eq order1}
X_{1:n}(t):=\max_{1\leq i\leq n}X_{i}(t)\geq \ldots \geq X_{m:n}(t)\geq \ldots \geq X_{n:n}(t):=\min_{1\leq i\leq n}X_{i}(t),\ \ t\ge 0.
\end{eqnarray}
The extremes properties of this process have attracted much attention in recent years.

The tail probability of $\sup_{t\in[0,T]} X_{m:n}(t)$ has been studied by  D\c{e}bicki, Hashorva, Ji and Ling (2015) and
 D\c{e}bicki, Hashorva, Ji and Tabi\'{s} (2014,2015).
 D\c{e}bicki et al. (2017)  obtained the limit distribution of $\sup_{t\in[0,T]} X_{m:n}(t)$.
\COM{Suppose that $r(t)$ satisfies the conditions (\ref{eq1.1}) and (\ref{eq1.2}). If $r\in[0,+\infty)$, we have
\begin{eqnarray}
\label{eq TB}
P\left\{a_{m,T}(\sup_{t\in[0,T]} X_{m:n}(t)-b_{m,T})\leq x\right\}\longrightarrow \mathbb{E}\exp\left(-e^{-x-r+\sqrt{2rm}\mathcal{N}}\right)
\end{eqnarray}
as $T\rightarrow\infty$ with $\mathcal{N}$ an $N(0,1)$ random variable, where the normalizing constants $a_{m,T}$ and $b_{m,T}$ are defined as follows
$$a_{m,T}=\sqrt{2m \ln T}, \quad b_{m,T}=\frac{1}{m}a_{m,T}+ a_{m,T}^{-1}\ln\Bigl(a_{m,T}^{2/\alpha-m}C_{n}^{m}\mathcal{H}_{m,\alpha}(2\pi)^{-m/2}\Bigr).$$
If $r=\infty$, and $\alpha\in(0,1]$, $r(t)$ is convex for all $t\geq 0$ with $\lim_{t\rightarrow\infty}r(t)=0$ and further $r(t)\ln t$ is monotone for large
$t$, then
\begin{eqnarray}
\label{eq TC}
P\left\{r^{-1/2}(T)(\sup_{t\in[0,T]} X_{m:n}(t)-\sqrt{1-r(T)}b_{m,T})\leq x\right\}\longrightarrow \Phi(x),
\end{eqnarray}
where $\Phi(\cdot)$ denotes the cumulative distribution function of a standard normal variable.}
For some related studies on extremes properties for Gaussian order statistics processes, we refere to D\c{e}bicki and Kosi\'{n}ski (2017) and Zhao (2017).

In applied fields, however, the above limit results
can not be used directly, since the available samples are usually
over a discrete set of times.
 Therefore, it is crucial
to investigate the asymptotic relation between the
maxima of the continuous time processes and the maxima of the processes sampled at discrete time points.

Piterbarg (2004) first studied the asymptotic relation between the
continuous time maximum and the discrete version maximum of stationary Gaussian processes. This type of results are called
Piterbarg's max-discretisation theorems in the literature, see e.g.  Tan and Hashorva (2014).
Piterbarg's max-discretisation theorems have been
extended to more general Gaussian cases, see H\"{u}sler (2004), H\"{u}sler and Piterbarg (2004), Tan and Hashorva (2014),
Hashorva and Tan (2015) and Tan and Wang (2015).
Although the Piterbarg's max-discretisation theorems for Gaussian processes
have been  studied extensively under different conditions in the past, it is far from complete.
Extending the above results to non-Gaussian case is also interesting, since most of reality can not
be modeled by Gaussian model.  Turkman (2012) considered this problem by adopting the model from Albin (1990),
but there are some mistakes in the paper, which have been corrected by Ling et al. (2017).
 Ling and Tan (2016) dealt with the problems for chi-processes.
The goal of this paper is to study the Piterbarg's max-discretisation theorems for order statistics processes.

Following Piterbarg (2004),
we consider uniform grids
$\mathfrak{R}(p)=\{kp: k\in \mathbb{N}\}$, $p=p(T)>0$.
A grid $\mathfrak{R}(p)$ is called sparse if $p$ is such that
$$p(T)(\frac{2}{m}\ln T)^{1/\alpha}\rightarrow D,\ \ T\rightarrow\infty$$
with $D=\infty$. If $D\in(0, \infty)$, the grid is a Pickands grid, and if $D = 0$, the grid is dense.
For any grid $\mathfrak{R}(p)$, define
$$M_{m}(T):= \sup_{t\in[0,T]} X_{m:n}(t) \ \ \mbox{and}\ \ M_{m}^{p}(T):=\sup_{kp\in[0,T]} X_{m:n}(kp).$$

The paper is organized as follows: In
Section 2, we present the main results. Section 3 gives the proofs. Some technical auxiliary results are presented in Section
4. Let  $\Psi(\cdot)$ and $\Phi(\cdot)$ denote the tail distribution function and cumulative distribution function of a standard normal
variable, respectively. 

\section{Main results}

For stating our main results,  we introduce the
following Pickands type constants. Let $B_{\alpha/2}^{(i)}$, $i=1,2,\ldots,m$ be independent fractional Brownian motions and define
$$\mathcal{H}_{m,\alpha}(\lambda)=\int_{\mathbb{R}^{m}}e^{\sum_{i=1}^{m}w_{i}}P\left\{\sup_{t\in[0,\lambda]}\min_{1\leq i\leq m}(\sqrt{2}B_{\alpha/2}^{(i)}(t)-t^{\alpha}>w_{i})\right\}d\mathbf{w}.$$
It follows from D\c{e}bicki, Hashorva, Ji and Tabi\'{s} (2015) that
$$\mathcal{H}_{m,\alpha}:=\lim_{\lambda\rightarrow\infty}\frac{\mathcal{H}_{m,\alpha}(\lambda)}{\lambda}\in (0,\infty).$$

\textbf{Theorem 2.1}. {\sl Let $\{X(t), t\geq0\}$ be a standard (zero-mean, unit-variance)
stationary Gaussian process with correlation functions $r(\cdot)$ satisfying (\ref{eq1.1}) and (\ref{eq1.2}) with $r\in[0,+\infty)$.
Then for any sparse grid $\mathfrak{R}(p)$,
\begin{eqnarray}
\label{eq2.1}
&&P\left\{a_{m,T}\big(M_{m}(T)-b_{m,T}\big)\leq x,
           a_{m,T}\big(M_{m}^{p}(T)- b_{m,T}^{p}\big)\leq y\right\}\nonumber \\
&&\ \ \ \ \ \ \longrightarrow \mathbb{E}\exp\left(-\big(e^{-x-r+\sqrt{2rm}\mathcal{N}}+e^{-y-r+\sqrt{2rm}\mathcal{N}}\big)\right)
\end{eqnarray}
as $T\rightarrow\infty,$ where $\mathcal{N}$ is an $N(0,1)$ random variable and the normalizing constants are defined as
$$a_{m,T}=\sqrt{2m \ln T}, \quad b_{m,T}=\frac{1}{m}a_{m,T}+ a_{m,T}^{-1}\ln\Bigl(a_{m,T}^{2/\alpha-m}C_{n}^{m}\mathcal{H}_{m,\alpha}(2\pi)^{-m/2}\Bigr)$$
and
$$b_{m,T}^{p}=\frac{1}{m}a_{m,T}+ a_{m,T}^{-1}\ln \Bigl(a_{m,T}^{-m} C_{n}^{m}(2\pi)^{-m/2}p^{-1}\Bigr)$$
with  $C_{n}^{m}=\frac{n!}{m!(n-m)!}$.
}

As a special case, we can obtain the limit distribution of the
maximum for discrete order statistics processes, which is of independent interest.

\textbf{Corollary 2.1}. {\sl Let $\{X(t), t\geq0\}$ be a standard (zero-mean, unit-variance)
stationary Gaussian process with correlation functions $r(\cdot)$ satisfying (\ref{eq1.1}) and (\ref{eq1.2}) with $r\in[0,+\infty)$.
Then for any sparse grid $\mathfrak{R}(p)$,
\begin{eqnarray}
\label{eq2.4}
&&P\bigg\{a_{m,T}\big(M_{m}^{p}(T)- b_{m,T}^{p}\big)\leq x\bigg\}
\longrightarrow \mathbb{E}\exp\left(-e^{-x-r+\sqrt{2rm}\mathcal{N}}\right)
\end{eqnarray}
as $T\rightarrow\infty.$ }

Before presenting the result for Pickands grids, we introduce the
following Pickands type constants.
For any $d>0$, define for $k\in \mathbb{N}$,
$$\mathcal{H}_{d,m,\alpha}(\lambda)=\int_{\mathbb{R}^{m}}e^{\sum_{i=1}^{m}w_{i}}P\left\{\sup_{kd\in[0,\lambda]}\min_{1\leq i\leq m}(\sqrt{2}B_{\alpha/2}^{(i)}(kd)-(kd)^{\alpha}>w_{i})\right\}d\mathbf{w}.$$
We have (by the same arguments as in  D\c{e}bicki, Hashorva, Ji and Tabi\'{s} (2015)),
$$\mathcal{H}_{d,m,\alpha}:=\lim_{\lambda\rightarrow\infty}\frac{\mathcal{H}_{d,m,\alpha}(\lambda)}{\lambda}\in(0,\infty).$$
For any $d>0$, define for $k\in \mathbb{N}$
\begin{eqnarray*}
&&\mathcal{H}_{d,m,\alpha}^{x,y}(\lambda)=\int_{\mathbb{R}^{m}}e^{\sum_{i=1}^{m}w_{i}}
P\bigg\{\sup_{t\in[0,\lambda]}\min_{1\leq i\leq m}(\sqrt{2}B_{\alpha/2}^{(i)}(t)-t^{\alpha}>w_{i}+x,\\
&&\ \ \ \ \ \ \ \ \ \ \ \ \ \ \ \ \ \ \ \ \ \ \ \ \ \ \ \ \ \ \sup_{kd\in[0,\lambda]}\min_{1\leq i\leq m}(\sqrt{2}B_{\alpha/2}^{(i)}(kd)-(kd)^{\alpha}>w_{i}+y\bigg\}d\mathbf{w}
\end{eqnarray*}
we have (see the proof in Appendix A)
$$\mathcal{H}_{d,m,\alpha}^{x,y}:=\lim_{\lambda\rightarrow\infty}
\frac{\mathcal{H}_{d,m,\alpha}^{x,y}(\lambda)}{\lambda}\in
(0,\infty).$$

\textbf{Theorem 2.2}. {\sl Let $\{X(t), t\geq0\}$ be a standard (zero-mean, unit-variance)
stationary Gaussian process with correlation functions $r(\cdot)$ satisfying (\ref{eq1.1}) and (\ref{eq1.2}) with $r\in[0,+\infty)$.
Then for any Pickands grid $\mathfrak{R}(p)=\mathfrak{R}(d(\frac{2}{m}\ln
T)^{-1/\alpha})$ with $d>0$,
\begin{eqnarray}
\label{eq2.2}
&&P\left\{a_{m,T}\big(M_{m}(T)-b_{m,T}\big)\leq x,
           a_{m,T}\big(M_{m}^{p}(T)- b_{d,m,T}\big)\leq y\right\}\nonumber\\
&&\ \ \ \ \ \longrightarrow \mathbb{E}\exp\left(-\big(e^{-x-r+\sqrt{2rm}\mathcal{N}}
+e^{-y-r+\sqrt{2rm}\mathcal{N}}-\mathcal{H}_{d,m,\alpha}^{\ln(\mathcal{H}_{m,\alpha})+x, \ln
(\mathcal{H}_{d,m,\alpha})+y}e^{-r+\sqrt{2rm}\mathcal{N}}\big)\right)
\end{eqnarray}
as $T\rightarrow\infty,$ where
$$b_{d,m,T}=\frac{1}{m}a_{m,T}+ a_{m,T}^{-1}\ln\Bigl(a_{m,T}^{2/\alpha-m}C_{n}^{m}\mathcal{H}_{d,m,\alpha}(2\pi)^{-m/2}\Bigr).$$ }

\textbf{Theorem 2.3}. {\sl Let $\{X(t), t\geq0\}$ be a standard (zero-mean, unit-variance)
stationary Gaussian process with correlation functions $r(\cdot)$ satisfying (\ref{eq1.1}) and (\ref{eq1.2}) with $r\in[0,+\infty)$.
 Then for any dense grid $\mathfrak{R}(p)$,
\begin{eqnarray}
\label{eq2.3}
&&P\left\{a_{m,T}\big(M_{m}(T)-b_{m,T}\big)\leq x, a_{m,T}\big(M_{m}^{p}(T)- b_{m,T}\big)\leq y\right\}
\longrightarrow \mathbb{E}\exp\left(-e^{-\min(x,y)-r+\sqrt{2rm}\mathcal{N}}\right)
\end{eqnarray}
as $T\rightarrow\infty.$ }

\textbf{Theorem 2.4}. {\sl Let $\{X(t), t\geq0\}$ be a standard (zero-mean, unit-variance)
stationary Gaussian process with correlation functions $r(\cdot)$ satisfying (\ref{eq1.1}) and (\ref{eq1.2}). Suppose that $r=\infty$
and $\alpha\in(0,1]$, $r(t)$ is convex for all $t\geq 0$ with $\lim_{t\rightarrow\infty}r(t)=0$ and further $r(t)\ln t$ is monotone for large
$t$. Then for any grid $\mathfrak{R}(p)$,
\begin{eqnarray}
\label{eq2.4}
&&P\left\{\frac{1}{\sqrt{r(T)}}(M_{m}(T)-\sqrt{1-r(T)}b_{m,T})\leq x, \frac{1}{\sqrt{r(T)}}(M_{m}^{p}(T)-\sqrt{1-r(T)}b_{m,T}^{*})\leq y\right\}
\longrightarrow \Phi(\min\{x,y\})
\end{eqnarray}
as $T\rightarrow\infty,$ where $b_{m,T}^{*}=b_{m,T}^{p}$ for a sparse grid; $b_{m,T}^{*}=b_{d,m,T}$ for a Pickands grid;
 $b_{m,T}^{*}=b_{m,T}$ for a dense grid.
}

\textbf{Remark 2.1}. We note that the above results under different grids still holds. It is not difficult to check
it by combining the method used in this paper with the one in Hashorva and Tan (2015).

\section{Proofs}
\subsection{Proofs of Theorems 2.1-2.3}

First, define $\rho(T)=r/\ln T$ and let $0<c<a<1$
be positive constants which will be determined by (\ref{bkk}) and (\ref{eqTB}) in Appendixes A and B, respectively.
Following Piterbarg (2004), divide $[0,T]$ into intervals with
length $T^{a}$ alternating with shorter intervals with length
$T^{c}$. Note that the numbers of the long intervals is at most
$l=l_{T}=\lfloor T/(T^{a}+T^{c})\rfloor$, where $\lfloor x\rfloor$ denotes the integral part of $x$.
Let
$\mathbf{O}_{i}=
[(i-1)(T^{a}+T^{c}),(i-1)(T^{a}+T^{c})+T^{a}]$,
 $\mathbf{Q}_{i}=
((i-1)(T^{a}+T^{c})+T^{a},i(T^{a}+T^{c}))$,
$i=1,\cdots,l$ and
$\mathbf{O}=\cup_{i} \mathbf{O}_{i}$.
There is still a remaining interval denoted by $\mathbf{Q}_{l+1}$ which  play no role in our consideration.

Let $\{Y^{i}_{j}(t), t\in \mathbf{E}_{i}\}$,
$i= 1,2, \cdots, l, j=1,2,\ldots,n$, be independent copies of
$\{X(t), t\geq 0\}$,
where $\mathbf{E}_{i}:=[(i-1)(T^{a}+T^{c}),i(T^{a}+T^{c}))$.
Further, define for $j=1,2,\ldots,n$
$$\xi^{T}_{j}(t)=\big(1-\rho(T)\big)^{1/2}\sum_{i=1}^{l}Y^{i}_{j}(t)\mathbb{I}(t\in \mathbf{E}_{i})+\rho^{1/2}(T)\mathcal{N},$$
where $\mathbb{I}$ is the indicator function and $\mathcal{N}$ is an $N(0,1)$ random variable, which is independent of $\{X(t), t\geq 0\}$ and $\{Y^{i}_{j}(t), t\in \mathbf{E}_{i}\}$, $i= 1,2, \cdots, l, j=1,2,\ldots,n$.
Denote by
$\varrho(t, s)$ the covariance function of
$\{\xi_{j}^{T}(t)\}$ and we have
\[
 \varrho(t, s)=\left\{
 \begin{array}{cc}
  {r(|t-s|)+(1-r(|t-s|))\rho(T)},    &s\in \mathbf{E}_{i}, t\in \mathbf{E}_{j}, i=j;\\
  {\rho(T)},    &s\in \mathbf{E}_{i}, t\in \mathbf{E}_{j}, i\neq j.
 \end{array}
  \right.
\]
Note that $\mathbf{O}_{i}\subset \mathbf{E}_{i}$, $i=1,2,\cdots,l$.
Let $Y_{m:n}^{i}(t)$ and $\xi^{T}_{m:n}(t)$ be the order statistics
processes formed by $Y_{j}^{i}(t)$ and $\xi^{T}_{j}(t)$, respectively.

In the sequel, $\mathcal{C}$ shall denote
positive constant whose values may vary from place to place.
For simplicity, define $$q=q(T)=b(\frac{2}{m}\ln T)^{-1/\alpha}$$
for some constant $b>0$.

We need the following lemmas to prove Theorems 2.1 and 2.2.

\textbf{Lemma 3.1}. {\sl Suppose that the grid
$\mathfrak{R}(p)$ is a sparse grid or Pickands grid. For any
$B>0$, we have for all $x,y\in[-B,B]$,
\begin{eqnarray*}
&&\bigg|P\left\{a_{m,T}\big(M_{m}(T)-b_{m,T}\big)\leq x,
           a_{m,T}\big(M_{m}^{p}(T)- b_{m,T}^{*}\big)\leq y\right\}\\
&&\ \ \ \ \ -P\left\{a_{m,T}\big(\max_{t\in \mathbf{O}}X_{m:n}(t)-b_{m,T}\big)\leq x,
a_{m,T}\big(\max_{t\in\mathfrak{R}(p)\cap\mathbf{O}}X_{m:n}(t)- b_{m,T}^{*}\big)\leq y\right\}\bigg|\rightarrow0
\end{eqnarray*}
as $T\rightarrow\infty$, where
$b_{m,T}^{*}=b_{m,T}^{p}$ for  sparse grid and
$b_{m,T}^{*}=b_{m,d,T}$ for  Pickands grid.
}

\textbf{Proof:} The proof is similar to that of Lemma 6 of Piterbarg
(2004). Clearly, we have
\begin{eqnarray}
\label{eq400}
&&\bigg|P\left\{a_{m,T}\big(M_{m}(T)-b_{m,T}\big)\leq x,
           a_{m,T}\big(M_{m}^{p}(T)- b_{m,T}^{*}\big)\leq y\right\}\nonumber\\
&&\ \ \ \ \ -P\left\{a_{m,T}\big(\max_{t\in \mathbf{O}}X_{m:n}(t)-b_{m,T}\big)\leq x,
a_{m,T}\big(\max_{t\in\mathfrak{R}(p)\cap\mathbf{O}}X_{m:n}(t)- b_{m,T}^{*}\big)\leq y\right\}\bigg|\nonumber\\
&&\ \ \ \ \ \leq \sum_{i=1}^{l+1} P\left\{\max_{t\in \mathbf{Q_{i}}}X_{m:n}(t)>b_{m,T}+x/a_{m,T}\right\}
+\sum_{i=1}^{l+1}P\left\{\max_{t\in\mathfrak{R}(p)\cap\mathbf{Q_{i}}}X_{m:n}(t)>b_{m,T}^{*}+y/a_{m,T}\right\}.
\end{eqnarray}
To bound the right hand side of (\ref{eq400}), we need the following result
for any $S\in(0,\exp(cu^{2}))$ for some $c\in(0,1/2)$
\begin{eqnarray}
\label{eq chi1}
P\left\{\sup_{t\in[0,S]} X_{m:n}(t) >u\right\}&=&C_{n}^{m}P\left\{\sup_{t\in[0,S]} X_{m:m}(t)>u\right\}(1+o(1))\nonumber\\
&=&C_{n}^{m}S\mathcal{H}_{m,\alpha}u^{\frac{2}{\alpha}}\Psi^{m}(u)(1+o(1)), \quad  u\rightarrow\infty.
\end{eqnarray}
For fixed $S$, it follows from D\c{e}bicki, Hashorva, Ji and Ling (2015). For the case $S\rightarrow\infty$,
it can be proved by the same arguments as in Lemma D.2 in Piterbarg (1996).
Thus, by the choice of $a_{m,T}$
and $b_{m,T}$, we have (denote by $mes(\cdot)$ the Lebesgue measure)
\begin{eqnarray*}
\sum_{i=1}^{l+1}P\left\{\max_{t\in \mathbf{Q_{i}}}X_{m:n}(t)>b_{m,T}+x/a_{m,T}\right\}
&=&O(1)\sum_{i=1}^{l+1}mes(\mathbf{O_{i}})(b_{m,T}+x/a_{m,T})^{2/\alpha}\Psi^{m}(b_{m,T}+x/a_{m,T})\\
&=&O(1)\frac{\sum_{i=1}^{l+1}mes(\mathbf{Q_{i}})}{T}\\
&\leq &O(1)\frac{(l+1)T^{c}}{T}\rightarrow0
\end{eqnarray*}
as $T\rightarrow\infty$. In light of  the second assertion in Lemmas A1 and A2 in the Appendix for a sparse grid and Pickands grid, respectively, we can
get the same estimation for the second probability in the right-hand side
of (\ref{eq400}), hence the proof is complete. \hfill $\Box$

\textbf{Lemma 3.2}. {\sl Suppose that the grids
$\mathfrak{R}(p)$ is a sparse grid or Pickands grid. For any
$B>0$, we have for all $x,y\in[-B,B]$ and the Pickands grids $\mathfrak{R}(q)=\mathfrak{R}(b(\frac{2}{m}\ln T)^{-1/\alpha})$
\begin{eqnarray}
\label{T1}
&&\bigg|P\left\{a_{m,T}\big(\max_{t\in \mathbf{O}}X_{m:n}(t)-b_{m,T}\big)\leq x,
a_{m,T}\big(\max_{t\in\mathfrak{R}(p)\cap\mathbf{O}}X_{m:n}(t)- b_{m,T}^{*}\big)\leq y\right\}\nonumber\\
&&\ \ \ \ \ -P\left\{a_{m,T}\big(\max_{t\in \mathfrak{R}(q)\cap\mathbf{O}}X_{m:n}(t)-b_{m,T}\big)\leq x,
a_{m,T}\big(\max_{t\in\mathfrak{R}(p)\cap\mathbf{O}}X_{m:n}(t)- b_{m,T}^{*}\big)\leq y\right\}\bigg|\rightarrow0
\end{eqnarray}
as $T\rightarrow\infty$ and $b\downarrow0$, where
$b_{m,T}^{*}=b_{m,T}^{p}$ for  sparse grids and
$b_{m,T}^{*}=b_{m,d,T}$ for  Pickands grids.
}

\textbf{Proof:} It is easy to see that the left hand side of (\ref{T1}) is bounded above by
\begin{eqnarray*}
&&P\left\{\max_{t\in \mathbf{O}}X_{m:n}(t)> a_{m,T}^{-1}x+b_{m,T},
\max_{t\in\mathfrak{R}(q)\cap\mathbf{O}}X_{m:n}(t)\leq a_{m,T}^{-1}x+b_{m,T}\right\}\\
&&\leq TP\left\{\max_{t\in [0,1]}X_{m:n}(t)> a_{m,T}^{-1}x+b_{m,T},
\max_{t\in\mathfrak{R}(q)\cap[0,1]}X_{m:n}(t)\leq a_{m,T}^{-1}x+b_{m,T}\right\}
\end{eqnarray*}
D\c{e}bicki, Hashorva, Ji and Ling (2017) has shown that the order statistics processes $X_{m:n}(t)$ satisfied Conditions $B$ and $C^{0}(\Lambda)$
in Albin (1990), which implies that condition (3.7) in Leadbetter and Rootz\'{e}n  (1982)  holds (see the proof of Theorem 10 in Albin (1990)).
Applying condition (3.7) in Leadbetter and Rootz\'{e}n (1982), we can see that the above probability does not exceed
$$o(1)T(b_{m,T}+x/a_{m,T})^{2/\alpha}\Psi^{m}(b_{m,T}+x/a_{m,T})=o(1)$$
as $T\rightarrow\infty$ and $b\downarrow0$ .\hfill$\Box$

\textbf{Lemma 3.3}. {\sl Suppose that the grid
$\mathfrak{R}(p)$ is a  sparse grid or  Pickands grid. For any
$B>0$ we have for all $x,y\in[-B,B]$ and the Pickands grids $\mathfrak{R}(q)=\mathfrak{R}(b(\frac{2}{m}\ln T)^{-1/\alpha})$
\begin{eqnarray*}
&&\bigg|P\left\{a_{m,T}\big(\max_{t\in \mathfrak{R}(q)\cap\mathbf{O}}X_{m:n}(t)-b_{m,T}\big)\leq x,
a_{m,T}\big(\max_{t\in\mathfrak{R}(p)\cap\mathbf{O}}X_{m:n}(t)- b_{m,T}^{*}\big)\leq y\right\}\\
&&\ \ \ \ \ -P\left\{a_{m,T}\big(\max_{t\in \mathfrak{R}(q)\cap\mathbf{O}}\xi_{m:n}^{T}(t)-b_{m,T}\big)\leq x,
a_{m,T}\big(\max_{t\in\mathfrak{R}(p)\cap\mathbf{O}}\xi_{m:n}^{T}(t)- b_{m,T}^{*}\big)\leq y\right\}\bigg|\rightarrow 0
\end{eqnarray*}
uniformly for $b>0$ as $T\rightarrow\infty$, where
$b_{m,T}^{*}=b_{m,T}^{p}$ for  sparse grids and
$b_{m,T}^{*}=b_{m,d,T}$ for  Pickands grids.
}

\textbf{Proof:} For the sake of simplicity, let $u_{T}=
b_{m,T}+x/a_{m,T}$, $u_{T}^{*}=
b_{m,T}^{*}+y/a_{m,T}$. Using Lemma C1 in Appendix C,
we have
\begin{eqnarray*}
\label{eq401}
&&\bigg|P\left\{a_{m,T}\big(\max_{t\in \mathfrak{R}(q)\cap\mathbf{O}}X_{m:n}(t)-b_{m,T}\big)\leq x,
a_{m,T}\big(\max_{t\in\mathfrak{R}(p)\cap\mathbf{O}}X_{m:n}(t)- b_{m,T}^{*}\big)\leq y\right\}\\
&&\ \ \ \ \ -P\left\{a_{m,T}\big(\max_{t\in \mathfrak{R}(q)\cap\mathbf{O}}\xi_{m:n}^{T}(t)-b_{m,T}\big)\leq x,
a_{m,T}\big(\max_{t\in\mathfrak{R}(p)\cap\mathbf{O}}\xi_{m:n}^{T}(t)- b_{m,T}^{*}\big)\leq y\right\}\bigg|\nonumber\\
&&\leq \mathcal{C}\sum_{t\in \mathfrak{R}(q)\cap\mathbf{O}_{i}, s\in \mathfrak{R}(q)\cap\mathbf{O}_{j},\atop t\neq s, 1\leq i,j\leq l }|r(|t-s|)-\varrho(s,t)|
\int_{0}^{1}\frac{u_{T}^{-2(m-1)}}{(1-r^{(h)}(t,s))^{m/2}}\exp\left(-\frac{mu_{T}^{2}}{1+|r^{(h)}(t,s)|}\right)dh\nonumber\\
&&+\mathcal{C}\sum_{t\in \mathfrak{R}(p)\cap\mathbf{O}_{i}, s\in \mathfrak{R}(p)\cap\mathbf{O}_{j},\atop t\neq s, 1\leq i,j\leq l }|r(|t-s|)-\varrho(s,t)|
\int_{0}^{1}\frac{u_{T}^{-2(m-1)}}{(1-r^{(h)}(t,s))^{m/2}}\exp\left(-\frac{m(u_{T}^{*})^{2}}{1+|r^{(h)}(t,s)|}\right)dh\nonumber\\
&&+\mathcal{C}\sum_{t\in \mathfrak{R}(q)\cap\mathbf{O}_{i}, s\in \mathfrak{R}(p)\cap\mathbf{O}_{j},\atop t\neq s, 1\leq i,j\leq l }|r(|t-s|)-\varrho(s,t)|
\int_{0}^{1}\frac{u_{T}^{-2(m-1)}}{(1-r^{(h)}(t,s))^{m/2}}\exp\left(-\frac{\frac{1}{2}m[u_{T}^{2}+(u_{T}^{*})^{2}]}{1+|r^{(h)}(t,s)|}\right)dh,
\end{eqnarray*}
where
$r^{(h)}(t,s)=hr(|t-s|)+(1-h)\varrho(t,s)$.
Now, the lemma follows from Lemmas B1-B3 in the Appendix B.
\hfill$\Box$

\textbf{Proof of Theorem 2.1}. First, by the definition of
$\{\xi_{m:n}^{T}(t)\}$, we have
\begin{eqnarray*}
\label{eq412}
&&P\left\{a_{m,T}\big(\max_{t\in \mathfrak{R}(q)\cap\mathbf{O}}\xi_{m:n}^{T}(t)-b_{m,T}\big)\leq x,
a_{m,T}\big(\max_{t\in\mathfrak{R}(p)\cap\mathbf{O}}\xi_{m:n}^{T}(t)- b_{m,T}^{*}\big)\leq y\right\}\nonumber\\
&&= \frac{1}{(2\pi)^{1/2}}\int_{\mathbb{R}}e^{-\frac{1}{2}z^{2}}\prod_{i=1}^{l}P\left\{\max_{t\in \mathfrak{R}(q)\cap\mathbf{O}_{i}}Y_{m:n}^{i}(t)\leq \frac{b_{m,T}+x/a_{m,T}-\rho^{1/2}(T)z}{(1-\rho(T))^{1/2}},\right.\\
  &&\ \ \ \left.   \max_{t\in \mathfrak{R}(p)\cap\mathbf{O}_{i}}Y_{m:n}^{i}(t)\leq \frac{b_{m,T}^{*}+x/a_{m,T}-\rho^{1/2}(T)z}{(1-\rho(T))^{1/2}}\right\}dz\\
&&=\mathbb{E}\prod_{i=1}^{l}
P\left\{\max_{t\in \mathfrak{R}(q)\cap\mathbf{O}_{i}}Y_{m:n}^{i}(t)\leq v_{T},
   \max_{t\in \mathfrak{R}(p)\cap\mathbf{O}_{i}}Y_{m:n}^{i}(t)\leq v_{T}^{*}\right\},
\end{eqnarray*}
 where
\begin{eqnarray}
\label{equ413}
v_{T}:=\frac{b_{m,T}+x/a_{m,T}-\rho^{1/2}(T)\mathcal{N}}{(1-\rho(T))^{1/2}}
=\frac{x+r-\sqrt{2rm}\mathcal{N}}{a_{m,T}}+b_{m,T}+o(a_{m,T}^{-1}),
\end{eqnarray}
and
\begin{eqnarray}
\label{equ414}
v_{T}^{*}:=\frac{b_{m,T}^{*}+y/a_{m,T}-\rho^{1/2}(T)\mathcal{N}}{(1-\rho(T))^{1/2}}
=\frac{y+r-\sqrt{2rm}\mathcal{N}}{a_{m,T}}+b_{m,T}^{*}+o(a_{m,T}^{-1})
\end{eqnarray}
as $T\rightarrow\infty$ with $b_{m,T}^{*}=b_{m,T}^{p}$ for  sparse grid and
$b_{m,T}^{*}=b_{m,d,T}$ for Pickands grid.

Now, from Lemmas 3.1-3.3, we know that in order to prove Theorem 2.1, it
suffices to show that
\begin{eqnarray*}
&&\bigg|\mathbb{E}\prod_{i=1}^{l}
P\left\{\max_{t\in \mathfrak{R}(q)\cap\mathbf{O}_{i}}Y_{m:n}^{i}(t)\leq v_{T},
   \max_{t\in \mathfrak{R}(p)\cap\mathbf{O}_{i}}Y_{m:n}^{i}(t)\leq v_{T}^{*}\right\}\\
&&\ \ \ \ \ \ \ -\mathbb{E}\exp\left(-\big(e^{-x-r+\sqrt{2rm}\mathcal{N}}+e^{-y-r+\sqrt{2rm}\mathcal{N}}\big)\right)\bigg|\rightarrow 0
\end{eqnarray*}
as $T\rightarrow\infty$, where $v_{T}$ and
$v_{T}^{*}$ are defined in (\ref{equ413}) and (\ref{equ414}), respectively.
Using the stationarity of $\{Y_{m:n}^{i}(t), t\in \mathbf{E}_{i}\}$ with respect to $t$, we have
\begin{eqnarray*}
&&\prod_{i=1}^{l}
P\left\{\max_{t\in \mathfrak{R}(q)\cap\mathbf{O}_{i}}Y_{m:n}^{i}(t)\leq v_{T},
   \max_{t\in \mathfrak{R}(p)\cap\mathbf{O}_{i}}Y_{m:n}^{i}(t)\leq v_{T}^{*}\right\}\\
&&=\left(P\left\{\max_{t\in \mathfrak{R}(q)\cap[0,T^{a}]}Y_{m:n}^{i}(t)\leq v_{T},
   \max_{t\in \mathfrak{R}(p)\cap[0,T^{a}]}Y_{m:n}^{i}(t)\leq v_{T}^{*}\right\}\right)^{l}\\
&&=\exp\left(l\ln\left(P\left\{\max_{t\in \mathfrak{R}(q)\cap[0,T^{a}]}Y_{m:n}^{i}(t)\leq v_{T},
   \max_{t\in \mathfrak{R}(p)\cap[0,T^{a}]}Y_{m:n}^{i}(t)\leq v_{T}^{*}\right\}\right)\right)\\
&&=\exp\left(-l\left(1-P\left\{\max_{t\in \mathfrak{R}(q)\cap[0,T^{a}]}Y_{m:n}^{i}(t)\leq v_{T},
   \max_{t\in \mathfrak{R}(p)\cap[0,T^{a}]}Y_{m:n}^{i}(t)\leq v_{T}^{*}\right\}\right)+R_{n}\right),
\end{eqnarray*}
where $R_{n}$ is the remainder of the Taylor expansion $\ln x=-(1-x+x^{2}+\cdots)$ for $0<x<1$.
Since by the definitions of $v_{T}$ and $v_{T}^{*}$
$$P_{n}:=P\left\{\max_{t\in \mathfrak{R}(q)\cap[0,T^{a}]}Y_{m:n}^{i}(t)\leq v_{T},
   \max_{t\in \mathfrak{R}(p)\cap[0,T^{a}]}Y_{m:n}^{i}(t)\leq v_{T}^{*}\right\}\rightarrow 1$$
as $T\rightarrow\infty$, we get that the remainder
$R_{n}$ can be estimated as
$R_{n}=o(n(1-P_{n}))$. Using Lemma A1 and letting $b\downarrow0$, we get that
\begin{eqnarray*}
&&l\left(1-P\left\{\max_{t\in \mathfrak{R}(q)\cap[0,T^{a}]}Y_{m:n}^{i}(t)\leq v_{T},
   \max_{t\in \mathfrak{R}(p)\cap[0,T^{a}]}Y_{m:n}^{i}(t)\leq v_{T}^{*}\right\}\right)\\
&&\thicksim lT^{a-1}\left(e^{-x-r+\sqrt{2rm}\mathcal{N}}+e^{-y-r+\sqrt{2rm}\mathcal{N}}\right)\\
&&\thicksim e^{-x-r+\sqrt{2rm}\mathcal{N}}+e^{-y-r+\sqrt{2rm}\mathcal{N}}
\end{eqnarray*}
as $T\rightarrow\infty$, which combined with the dominated convergence theorem completes the
proof of Theorem 2.1.\hfill$\Box$

\textbf{Proof of Theorem 2.2}. As for the proof of Theorem 2.1, in view of Lemmas 3.1-3.3 in order to establish the proof we only
need to show
\begin{eqnarray*}
&&\bigg|\mathbb{E}\prod_{i=1}^{l}
P\left\{\max_{t\in \mathfrak{R}(q)\cap\mathbf{O}_{i}}Y_{m:n}^{i}(t)\leq v_{T},
   \max_{t\in \mathfrak{R}(p)\cap\mathbf{O}_{i}}Y_{m:n}^{i}(t)\leq v_{T}^{*}\right\}\\
&&\ \ \ \ \ \ \ -\mathbb{E}\exp\left(-\big(e^{-x-r+\sqrt{2rm}\mathcal{N}}+e^{-y-r+\sqrt{2rm}\mathcal{N}}
       -\mathcal{H}_{d,m,\alpha}^{\ln \mathcal{H}_{m,\alpha}+x,\ln \mathcal{H}_{d,m,\alpha}+y}e^{-r+\sqrt{2rm}\mathcal{N}}\big)\right)\bigg|\rightarrow 0
\end{eqnarray*}
as $T\rightarrow\infty$, where $v_{T}$ and
$v_{T}^{*}$ are defined in (\ref{equ413}) and (\ref{equ414}), respectively. Similar to the proof
of Theorem 2.1, using Lemma A2 and letting $b\downarrow0$, we get
\begin{eqnarray*}
&&l\left(1-P\left\{\max_{t\in \mathfrak{R}(q)\cap[0,T^{a}]}Y_{m:n}^{i}(t)\leq v_{T},
   \max_{t\in \mathfrak{R}(p)\cap[0,T^{a}]}Y_{m:n}^{i}(t)\leq v_{T}^{*}\right\}\right)\\
&&=lP\left\{\max_{t\in \mathfrak{R}(q)\cap[0,T^{a}]}Y_{m:n}^{i}(t)> v_{T}\right\}
   +lP\left\{\max_{t\in \mathfrak{R}(p)\cap[0,T^{a}]}Y_{m:n}^{i}(t)> v_{T}^{*}\right\}\\
&&\ \ -lP\left\{\max_{t\in \mathfrak{R}(q)\cap[0,T^{a}]}Y_{m:n}^{i}(t)> v_{T},
   \max_{t\in \mathfrak{R}(p)\cap[0,T^{a}]}Y_{m:n}^{i}(t)> v_{T}^{*}\right\}\\
&&\thicksim lT^{a-1}\left(e^{-x-r+\sqrt{2rm}\mathcal{N}}+e^{-y-r+\sqrt{2rm}\mathcal{N}}\right)\\
&&\ \ -lT^{a-1}\mathcal{H}_{d,m,\alpha}^{\ln
(\mathcal{H}_{m,\alpha})+x, \ln(\mathcal{H}_{d,m,\alpha})+y}e^{-r+\sqrt{2rm}\mathcal{N}}\\
&&\thicksim e^{-x-r+\sqrt{2rm}\mathcal{N}}+e^{-y-r+\sqrt{2rm}\mathcal{N}}-\mathcal{H}_{d,m,\alpha}^{\ln
(\mathcal{H}_{m,\alpha})+x, \ln(\mathcal{H}_{d,m,\alpha})+y}e^{-r+\sqrt{2rm}\mathcal{N}}
\end{eqnarray*}
as $T\rightarrow\infty$.
This and the dominated convergence theorem conclude the proof of
Theorem 2.2.\hfill$\Box$

\textbf{Proof of Theorem 2.3}. The proof  is the same as that of Theorem 2.2(ii) of Hashorva and Tan (2015), so we omit it. \hfill$\Box$

\subsection{Proof of Theorem 2.4}

The following lemma is needed for the proof of Theorem 2.4.  Note that by Polya's criterion, the convexity of $r(t)$ ensures that there is a separable stationary Gaussian process $\{Z^{T}(t), 0\leq t\leq T\}$ with
correlation function
$$\gamma_{T}(t)=(r(t)-r(T))/(1-r(T))\ \ \ \mbox{for}\ \ \ t\leq T.$$
Let $\{Z^{T}_{j}(t), 0\leq t\leq T\}$, $j=1,2,\ldots,n$, be independent copies of $Z^{T}(t)$
and $\{Z^{T}_{m:n}(t), 0\leq t\leq T\}$ be the  corresponding order statistics processes.
Let $$M_{m}^{Z}(T)=\max_{0\leq t\leq T}Z^{T}_{m:n}(t),\ \ M_{m}^{Z,p}(T)= \max_{t\in \mathfrak{R}(p)\cap[0,T]}Z^{T}_{m:n}(t).$$

\textbf{Lemma 3.4}. {\sl Let $Z^{T}(t)$ be defined as before. Under the conditions of Theorem 2.4, for any $\varepsilon>0$, as $T\rightarrow\infty$, we have
\begin{eqnarray}
\label{eqL3.1}
P\left\{|M_{m}^{Z}(T)-b_{m,T}|>\varepsilon r^{1/2}(T)\right\}\rightarrow 0,
\end{eqnarray}
\begin{eqnarray}
\label{eqL3.2}
P\left\{|M_{m}^{Z,p}(T)-b_{m,T}^{p}|>\varepsilon r^{1/2}(T)\right\}\rightarrow 0
\end{eqnarray}
for a sparse grid, and
\begin{eqnarray}
\label{eqL3.3}
P\left\{|M_{m}^{Z,p}(T)-b_{d,m,T}|>\varepsilon r^{1/2}(T)\right\}\rightarrow 0
\end{eqnarray}
for a Pickands grid, where $b_{m,T}, b_{m,T}^{p}$ and $b_{d,m,T}$ are defined as before.
}

\textbf{Proof:} We first show (\ref{eqL3.1}).
Note that $\gamma_{T}(t)$ satisfies,
\begin{eqnarray}
\label{eqL3.4}
\gamma_{T}(t)=\frac{r(t)-r(T)}{1-r(T)}=1-C(T)|t|^{\alpha}+o(|t|^{\alpha})\ \ \mbox{as}\ \  t\rightarrow 0,
\end{eqnarray}
where
$$C(T)=\frac{1}{1-r(T)}\rightarrow 1\ \ \mbox{as}\ \  T\rightarrow \infty.$$
Using the stationarity of $\{Z^{T}(t), 0\leq t\leq T\}$, (\ref{eq chi1}) and the definition of $b_{m,T}$, we have
\begin{eqnarray}
\label{eqL3.5}
P\left\{M_{m}^{Z}(T)-b_{m,T}>\varepsilon r^{1/2}(T)\right\}
&\leq & (\lfloor T\rfloor+1) P\left\{\max_{0\leq t\leq 1}Z^{T}_{m:n}(t)>\varepsilon r^{1/2}(T)+b_{m,T}\right\}\nonumber\\
&\leq&O(1)  (\lfloor T\rfloor+1)  (\varepsilon r^{1/2}(T)+b_{m,T})^{\frac{2}{\alpha}-m}e^{-\frac{1}{2}m(r^{1/2}(T)+b_{m,T})^{2}}\nonumber\\
&\leq& O(1)(\lfloor T\rfloor+1)  (\ln T)^{1/\alpha-m/2}e^{-\frac{1}{2}m(\frac{2}{m}\ln T +\frac{2}{m}(1/\alpha-m/2)\ln\ln T+\frac{2\sqrt{2}}{\sqrt{m}}(r(T)\ln T)^{1/2})}\nonumber\\
&=&O(1) e^{-\sqrt{2m}(r(T)\ln T)^{1/2}}.
\end{eqnarray}
Now, by the condition that $r(T)\ln T\uparrow \infty$, we get that
$$P\left\{M_{m}^{Z}(T)-b_{m,T}>\varepsilon r^{1/2}(T)\right\}\rightarrow 0\ \ \mbox{as}\ \  T\rightarrow \infty.$$
On the other hand, the following asymptotic relation has been proved in the proof of Theorem 4.1 (b) of D\c{e}bicki, Hashorva, Ji and Ling (2017),
$$P\left\{M_{m}^{Z}(T)-b_{m,T}<-\varepsilon r^{1/2}(T)\right\}\rightarrow 0\ \ \mbox{as}\ \  T\rightarrow \infty.$$
The proof of (\ref{eqL3.1}) is complete.

Next, we show (\ref{eqL3.2}). Let $\{\eta(t),t\geq 0\}$ be a standardized Gaussian process with covariance function $\varrho$ and
$\{\eta_{j}(t), t\geq 0\}$, $j=1,2,\ldots,n$ be independent copies of $\{\eta(t),t\geq 0\}$.
Let $\eta_{m:n}(t)$ be the order statistics process formed by $\{\eta_{j}(t), t\geq 0\}$, $j=1,2,\ldots,n$ and define
$$M_{m}^{\eta,p}(T,\varrho)=\max_{t\in\mathfrak{R}(p)\cap[0,T]}\eta_{m:n}(t).$$
Since $\gamma_{T}(t)>0$, by Normal comparison Lemma for order statistics (see Corollary 2.3 of D\c{e}bicki, Hashorva, Ji and Ling (2017)), we get
$$P\left\{M_{m}^{Z,p}(T)-b_{m,T}^{p}>\varepsilon r^{1/2}(T)\right\}\leq P\left\{M_{m}^{\eta,p}(T,0)-b_{m,T}^{p}>\varepsilon r^{1/2}(T)\right\}$$
Note that, by the definitions, $b_{m,T}^{p}r^{1/2}(T)\rightarrow\infty$ as $T\rightarrow\infty$. It follows from
Corollary 2.1 that
$$P\left\{M_{m}^{\eta,p}(T,0)-b_{m,T}^{p}>\varepsilon r^{1/2}(T)\right\}\rightarrow 0\ \ \mbox{as}\ \  T\rightarrow \infty,$$
Thus
$$P\left\{M_{m}^{Z,p}(T)-b_{m,T}^{p}>\varepsilon r^{1/2}(T)\right\}\rightarrow 0\ \ \mbox{as}\ \  T\rightarrow \infty.$$
By the same arguments as in the proof of Theorem 4.1 (b) of D\c{e}bicki, Hashorva, Ji and Ling (2017), we can show that
$$P\left\{M_{m}^{Z,p}(T)-b_{m,T}^{p}<-\varepsilon r^{1/2}(T)\right\}\rightarrow 0\ \ \mbox{as}\ \  T\rightarrow \infty,$$
which completes the proof of (\ref{eqL3.2}).  The proof of (\ref{eqL3.3}) is similar to that of (\ref{eqL3.1}), so we omit the details.

\textbf{Proof of Theorem 2.4}. For the cases of sparse and Pickands grid, as in the paper of
Mittal and Ylvisaker (1975),
represent $M_{m}(T)$ and $M_{m}^{p}(T)$ by
$$M_{m}(T)=(1-r(T))^{1/2}M_{m}^{Z}(T)+r^{1/2}(T)\mathcal{N}$$
and
$$M_{m}^{p}(T)=(1-r(T))^{1/2}M_{m}^{Z,p}(T)+r^{1/2}(T)\mathcal{N},$$
where $\mathcal{N}$ is a standard normal variable independent of $Z_{m:n}(t)$.
Using Lemma 3.4, we get that
\begin{eqnarray*}
&&P\left\{\frac{M_{m}(T)-(1-r(T))^{1/2}b_{m,T}}{r^{1/2}(T)}\leq x, \frac{M_{m}^{p}(T)-(1-r(T))^{1/2}b_{m,T}^{*}}{r^{1/2}(T)}\leq y\right\}\\
&&=P\left\{\frac{(1-r(T))^{1/2}(M_{m}^{Z}(T)-b_{T})}{r^{1/2}(T)}+\mathcal{N}\leq x, \frac{(1-r(T))^{1/2}(M_{m}^{Z,p}(T)-b_{T}^{*})}{r^{1/2}(T)}+\mathcal{N}\leq y\right\}\\
&&\rightarrow P\left\{\mathcal{N}\leq x, \mathcal{N}\leq y\right\}\\
&&=\Phi\left(\min\{x,y\}\right)
\end{eqnarray*}
as $T\rightarrow\infty$. The proof for the dense grid  is same as that of Theorem 2.2(ii) of Hashorva and Tan (2015), so we omit it.

\section{Appendix}
\subsection{Appendix A}

In this subsection, we give two auxiliary lemmas, which are used in the proofs of Theorem 2.1 and 2.2, respectively.
Let $Y_{j}(t)$, $j=1,2,\ldots, n$ be independent copies of $X(t)$ and $Y_{m:n}(t)$ be the corresponding order statistics processes.
The following fact will be extensively used in the proof. From assumption (\ref{eq1.1}), we can choose an $\epsilon>0$ such that for all
$|s-t|\leq\epsilon<2^{-1/\alpha}$
\begin{eqnarray}
\label{eq.A1}
\frac{1}{2}|s-t|^{\alpha}\leq 1-r(|t-s|)\leq 2|s-t|^{\alpha}.
\end{eqnarray}
Now, let $\vartheta(x)=\sup_{x\leq |t-s|\leq T}r(|t-s|)$.
Assumption (\ref{eq1.1}) implies that
$\vartheta(\epsilon)<1$ for all $T$ and any
$\epsilon\in(0,2^{-1/\alpha})$. Consequently, we may choose some
positive constants $a, c$ such that
\begin{eqnarray}
\label{bkk}
0 < c<a<\frac{1-\vartheta(\epsilon)}{1+\vartheta(\epsilon)} < 1
\end{eqnarray}
for all sufficiently large $T$.\\

\textbf{Lemma A1}. {\sl Let $q=q(T)=b(\frac{2}{m}\ln T)^{-1/\alpha}$
for some constant $b>0$. Under the conditions of Theorem 2.1
we have
$$P\left\{\max_{t\in \mathfrak{R}(q)\cap[0,T^{a}]}Y_{m:n}(t)> v_{T}\right\}=T^{a-1}e^{-x-r+\sqrt{2rm}\mathcal{N}}(1+o(1)),$$
as $b\downarrow0$ and $T\rightarrow\infty$;
$$P\left\{\max_{t\in\mathfrak{R}(p)\cap [0,T^{a}]}Y_{m:n}(t)> v_{T}^{*}\right\}
=T^{a-1}e^{-y-r+\sqrt{2rm}\mathcal{N}}(1+o(1)),$$
as  $T\rightarrow\infty$ and
$$P\left\{\max_{t\in \mathfrak{R}(q)\cap[0,T^{a}]}Y_{m:n}(t)> v_{T},
\max_{t\in\mathfrak{R}(p)\cap [0,T^{a}]}Y_{m:n}(t)> v_{T}^{*}\right\}= o(T^{a-1})$$
uniformly for $b>0$ as $T\rightarrow\infty$, where  $v_{T}$ and $v_{T}^{*}$ are defined in (\ref{equ413}) and (\ref{equ414}), respectively. }

\textbf{Proof:}
For the first assertion, recalling that $\mathfrak{R}(q)$ is a Pickands grid with $q=q(T)=b(\frac{2}{m}\ln T)^{-1/\alpha}$ and noting that
\begin{eqnarray*}
P\left\{\max_{t\in \mathfrak{R}(q)\cap[0,T^{a}]}Y_{m:n}(t)> v_{T}\right\}
&\rightarrow &P\left\{\max_{t\in [0,T^{a}]}Y_{m:n}(t)> v_{T}\right\}
\end{eqnarray*}
as $b\downarrow0$, the result follows from (\ref{eq chi1}) and definition of $v_{T}$ by some simple computations.
Next, we show the second assertion. Note that $\mathfrak{R}(p)$ is a sparse grid in this case.
By  Bonferroni inequality for large $T$
\begin{eqnarray*}
\sum_{t\in\mathfrak{R}(p)\cap[0,T^{a}]}P\left\{Y_{m:n}(t)> v_{T}^{*}\right\}
&\geq &P\left\{\max_{t\in\mathfrak{R}(p)\cap [0,T^{a}]}Y_{m:n}(t)> v_{T}^{*}\right\}\\
&\geq &\sum_{t\in\mathfrak{R}(p)\cap[0,T^{a}]}P\left\{Y_{m:n}(t)> v_{T}^{*}\right\}\\
&&\ \ \ -\sum_{t\in\mathfrak{R}(p)\cap [0,T^{a}]\atop
s\in\mathfrak{R}(p)\cap [0,T^{a}]}P\left\{Y_{m:n}(t)> v_{T}^{*},Y_{m:n}(s)> v_{T}^{*}\right\}\\
&=:& P_{T,1}-P_{T,2}.
\end{eqnarray*}
In view of Lemma 1 of  D\c{e}bicki, Hashorva, Ji and Ling (2015), we have
\begin{eqnarray*}
P\left\{Y_{m:n}(0)> u\right\}&=&C_{n}^{m}(P\left\{Y(0)> u\right\})^{m}(1+o(1))\\
&=&C_{n}^{m}\Psi^{m}(u)(1+o(1)),
\end{eqnarray*}
as $u\rightarrow\infty$. Therefore, by the definition of $v_{T}^{*}$, we have
\begin{eqnarray*}
P_{T,1}&=&\sum_{t\in\mathfrak{R}(p)\cap[0,T^{a}]}P\left\{Y_{m:n}(t)> v_{T}^{*}\right\}\\
&=&T^{a}p^{-1}P\left\{Y_{m:n}(0)> v_{T}^{*}\right\}\\
&=&T^{a}p^{-1}C_{n}^{m}\Psi^{m}(v_{T}^{*})(1+o(1))\\
&=&T^{a-1}e^{-y-r+\sqrt{2rm}\mathcal{N}}(1+o(1)),
\end{eqnarray*}
as $T\rightarrow\infty$, whereas to complete the proof, we only need to show $P_{T,2}=o(T^{a-1})$ uniformly for $b>0$ as $T\rightarrow\infty$.
Split the term $P_{T,2}$ into two parts as
\begin{eqnarray*}
P_{T,2}&=&
\sum_{t\in\mathfrak{R}(p)\cap [0,T^{a}],
s\in\mathfrak{R}(p)\cap [0,T^{a}]\atop |t-s|\leq \epsilon}P\left\{Y_{m:n}(t)> v_{T}^{*},Y_{m:n}(s)> v_{T}^{*}\right\}\\
&&+\sum_{t\in\mathfrak{R}(p)\cap [0,T^{a}],
s\in\mathfrak{R}(p)\cap [0,T^{a}] \atop |t-s|> \epsilon}P\left\{Y_{m:n}(t)> v_{T}^{*},Y_{m:n}(s)> v_{T}^{*}\right\}\\
&=:& P_{T,21}+P_{T,22},
\end{eqnarray*}
where $\epsilon$ is chosen such that (\ref{eq.A1}) holds.
By Lemma 1 of  D\c{e}bicki, Hashorva, Ji and Ling (2015) and the proof of Lemma 9 in D\c{e}bicki, Hashorva, Ji and Tabi\'{s} (2015), we have
\begin{eqnarray*}
\label{eq.A12}
P_{T,21}&=&\sum_{t\in\mathfrak{R}(p)\cap [0,T^{a}],
s\in\mathfrak{R}(p)\cap [0,T^{a}]\atop |t-s|\leq \epsilon}P\left\{Y_{m:n}(s)> v_{T}^{*}\right\}P\left\{Y_{m:n}(t)> v_{T}^{*}|Y_{m:n}(s)> v_{T}^{*}\right\}\\
&\leq& \sum_{t\in\mathfrak{R}(p)\cap [0,T^{a}],
s\in\mathfrak{R}(p)\cap [0,T^{a}]\atop |t-s|\leq \epsilon}C_{n}^{m}\Psi^{m}(v_{T}^{*})2^{m+1}\Psi^{m}\left(v_{T}^{*}\sqrt{\frac{1-r(s-t)}{1+r(s-t)}}\right)
\end{eqnarray*}
By (\ref{eq.A1}), we can choose $\epsilon>0$ small enough such that
$$\frac{1-r(s-t)}{1+r(s-t)}\geq \frac{1}{4}|t-s|^{\alpha}$$
and we thus have
\begin{eqnarray*}
\label{eq.A12}
P_{T,21}&\leq&\mathcal{C}\sum_{t\in\mathfrak{R}(p)\cap [0,T^{a}],
s\in\mathfrak{R}(p)\cap [0,T^{a}]\atop |t-s|\leq \epsilon}\left[\Psi(v_{T}^{*})
\Psi\left(v_{T}^{*}\frac{1}{2}|t-s|^{\alpha/2}\right)\right]^{m}\\
&\leq& \mathcal{C}\Psi^{m}(v_{T}^{*})\sum_{t\in\mathfrak{R}(p)\cap [0,T^{a}],
s\in\mathfrak{R}(p)\cap [0,T^{a}]\atop |t-s|\leq \epsilon}
\frac{1}{|t-s|^{m\alpha/2}(v_{T}^{*})^{m}}\exp\left(-\frac{1}{8}m|t-s|^{\alpha}(v_{T}^{*})^{2}\right).
\end{eqnarray*}
Using  the definition of $v_{T}^{*}$ we obtain
\begin{eqnarray*}
P_{T,21}&\leq & \mathcal{C}T^{a}p^{-1}\Psi^{m}(v_{T}^{*})\sum_{0< kp\leq \epsilon}\frac{1}{(kp)^{m\alpha/2}(v_{T}^{*})^{m}}\exp\left(-\frac{1}{8}m(kp)^{\alpha}(v_{T}^{*})^{2}\right)\\
&=&\mathcal{C} T^{a-1}\sum_{0< kp\leq \epsilon}\frac{1}{[kp(\ln T)^{1/\alpha}]^{m\alpha/2}}\exp\left(-\frac{1}{4}m[kp(\ln T)^{1/\alpha}]^{\alpha}\right)(1+o(1))\\
&\leq &\mathcal{C} T^{a-1}\frac{1}{[p(\ln T)^{1/\alpha}]^{m\alpha/2}}\sum_{0<k\leq \lfloor\epsilon/p\rfloor+1}\exp\left(-\frac{1}{4}m[kp(\ln T)^{1/\alpha}]^{\alpha}\right)(1+o(1))\\
&\leq &\mathcal{C} T^{a-1}\frac{1}{[p(\ln T)^{1/\alpha}]^{m\alpha/2}}(1+o(1))\\
&=&T^{a-1}o(1),
\end{eqnarray*}
where we used additionally  the fact that $\lim_{T\rightarrow\infty}(\ln T)^{1/\alpha}p=\infty$, since $\mathfrak{R}(p)$ is a sparse grid. Thus,
we have $P_{T,21}=o(T^{a-1})$ as $T\rightarrow\infty$.

For the second term, by the Comparison Lemma for order statistics (see Theorem 2.4 in D\c{e}bicki, Hashorva, Ji and Ling (2017)), we have
\begin{eqnarray*}
P_{T,22}&=& \sum_{t\in\mathfrak{R}(p)\cap [0,T^{a}],
s\in\mathfrak{R}(p)\cap [0,T^{a}] \atop |t-s|> \epsilon}P\left\{Y_{m:n}(t)> v_{T}^{*},Y_{m:n}(s)> v_{T}^{*}\right\}\\
&\leq &  \sum_{t\in\mathfrak{R}(p)\cap [0,T^{a}],
s\in\mathfrak{R}(p)\cap [0,T^{a}] \atop |t-s|> \epsilon}\bigg[\Psi^{2m}(v_{T}^{*})
+\mathcal{C}(v_{T}^{*})^{-2(m-1)}\exp\left(-\frac{m(v_{T}^{*})^{2}}{1+|r(|t-s|)|}\right)\bigg]\\
&\leq &\mathcal{C}T^{a}p^{-1} \sum_{\epsilon\leq kp\leq T^{a}}\bigg[\Psi^{2m}(v_{T}^{*})
+\mathcal{C}(v_{T}^{*})^{-2(m-1)}\exp\left(-\frac{m(v_{T}^{*})^{2}}{1+|r(kp)|}\right)\bigg]\\
&\leq & \mathcal{C}T^{2a}p^{-2}\bigg[\Psi^{2m}(v_{T}^{*})
+\mathcal{C}(v_{T}^{*})^{-2(m-1)}\exp\left(-\frac{m(v_{T}^{*})^{2}}{1+\vartheta(\epsilon)}\right)\bigg]\\
&=:&P_{T,221}+P_{T,222}.
\end{eqnarray*}
Utilising again the fact that  $v_{T}^{*}\sim u_{T}\sim (\frac{2}{m}\ln T)^{1/2}$, we have
\begin{eqnarray*}
P_{T,221}&\leq& \mathcal{C}T^{2a}p^{-2}u_{T}^{-2m}\left[\exp\left(-\frac{1}{2}mu_{T}^{2}\right)\right]^{2}\\
&\leq& \mathcal{C}T^{2a}p^{-2}u_{T}^{-2m}T^{-2}\\
&=&o(T^{a-1})
\end{eqnarray*}
and
\begin{eqnarray*}
\label{eq.A16}
P_{T,222}&\leq&\mathcal{C}T^{2a}p^{-2}u_{T}^{-2(m-1)}\exp\left(-\frac{m(v_{T}^{*})^{2}}{1+\vartheta(\epsilon)}\right)\\
&\leq & \mathcal{C}T^{2a}p^{-2}u_{T}^{-2(m-1)}T^{-\frac{2}{(1+\vartheta(\epsilon))}}\nonumber\\
&\leq & \mathcal{C} T^{a-1} T^{a-\frac{1-\vartheta(\epsilon)}{1+\vartheta(\epsilon)}}p^{-2}(\ln T)^{-(m-1)}.
\end{eqnarray*}
Both (\ref{bkk}) and  $(\ln T)^{1/\alpha}p=\infty$ imply
$ S_{T,22}=o(T^{a-1})$ as $T\rightarrow\infty$. This completes the proof of the second assertion.

Now, we prove the third assertion.
Obviously, we have
\begin{eqnarray*}
&&P\left\{\max_{t\in \mathcal{R}(q)\cap[0,T^{a}]}Y_{m:n}(t)> v_{T},
\max_{t\in\mathfrak{R}(p)\cap [0,T^{a}]}Y_{m:n}(t)> v_{T}^{*}\right\}\\
&&\leq\sum_{t\in\mathfrak{R}(q)\cap [0,T^{a}],
s\in\mathfrak{R}(p)\cap [0,T^{a}],\atop |t-s|\leq\epsilon}\mathbb{P}\left\{Y_{m:n}(t)> v_{T},Y_{m:n}(s)> v_{T}^{*}\right\}\\
&&+\sum_{t\in\mathfrak{R}(q)\cap [0,T^{a}],
s\in\mathfrak{R}(p)\cap [0,T^{a}],\atop |t-s|>\epsilon}\mathbb{P}\left\{Y_{m:n}(t)> v_{T},Y_{m:n}(s)> v_{T}^{*}\right\}\\
&&=: Q_{T,21}+Q_{T,22}.
\end{eqnarray*}
By the same argument as for the term $P_{T,21}$, we have for
\begin{eqnarray*}
\label{eq.A12}
Q_{T,21}&\leq&\mathcal{C}\sum_{t\in\mathfrak{R}(q)\cap [0,T^{a}],
s\in\mathfrak{R}(p)\cap [0,T^{a}],\atop |t-s|\leq\epsilon}\left[\Psi^{m}(v_{T})
\Psi^{m}\left(v_{T}^{*}(\frac{1}{4}|t-s|^{\alpha}\right)\right]\\
&\leq& \mathcal{C}\Psi^{m}(v_{T})\sum_{t\in\mathfrak{R}(q)\cap [0,T^{a}],
s\in\mathfrak{R}(p)\cap [0,T^{a}],\atop |t-s|\leq\epsilon}
\frac{1}{|t-s|^{m\alpha/2}(v_{T}^{*})^{m}}\exp\left(-\frac{1}{8}m|t-s|^{\alpha}(v_{T}^{*})^{2}\right)\\
&\leq & \mathcal{C}T^{a}b^{-1}u_{T}^{2/\alpha}\Psi^{m}(v_{T})\sum_{0< kp\leq \epsilon}\frac{1}{(kp)^{m\alpha/2}(v_{T}^{*})^{m}}\exp\left(-\frac{1}{8}m(kp)^{\alpha}(v_{T}^{*})^{2}\right)\\
&\leq&\mathcal{C} T^{a-1}b^{-1}\sum_{0< kp\leq \epsilon}\frac{1}{(kp)^{m\alpha/2}(\ln T)^{m/2}}\exp\left(-\frac{1}{4}(kp)^{\alpha}\ln T\right)\\
&\leq &\mathcal{C} T^{a-1}b^{-1}\frac{1}{[(\ln T)^{1/2}p^{\alpha/2}]^{m}}\sum_{0<k\leq \lfloor\epsilon/p\rfloor+1}
\exp\left(-\frac{1}{4}(kp)^{\alpha}\ln T\right)\\
&\leq &\mathcal{C} T^{a-1}b^{-1}\frac{1}{[p(\ln T)^{1/\alpha}]^{m\alpha/2}}\\
&=&T^{a-1}o(1),
\end{eqnarray*}
uniformly for $b>0$, where we used additionally  the fact that $\lim_{T\rightarrow\infty}p(\ln T)^{1/\alpha}=\infty$,
since $\mathfrak{R}(p)$ is a  sparse grid.

To bound the term $Q_{T,22}$, by the Comparison Lemma for order statistics again, with the same arguments as for the term $P_{T,22}$,
 we have
\begin{eqnarray*}
Q_{T,22}&=& \sum_{t\in\mathfrak{R}(q)\cap [0,T^{a}],
s\in\mathfrak{R}(p)\cap [0,T^{a}],\atop |t-s|>\epsilon}\mathbb{P}\left\{Y_{m:n}(t)> v_{T},Y_{m:n}(s)> v_{T}^{*}\right\}\\
&\leq& \sum_{t\in\mathfrak{R}(q)\cap [0,T^{a}],
s\in\mathfrak{R}(p)\cap [0,T^{a}],\atop |t-s|>\epsilon}\mathbb{P}\left\{Y_{m:n}(t)> v_{T}^{*},Y_{m:n}(s)> v_{T}^{*}\right\}\\
&\leq& \mathcal{C}T^{2a}p^{-1}q^{-1}\bigg[\Psi^{2m}(v_{T}^{*})
+\mathcal{C}(v_{T}^{*})^{-2(m-1)}\exp\left(-\frac{m(v_{T}^{*})^{2}}{1+\vartheta(\epsilon)}\right)\bigg]\\
&=:&Q_{T,221}+Q_{T,222}.
\end{eqnarray*}
By the same arguments as for $P_{T,221}$ and $P_{T,222}$,
we can show that $Q_{T,221}=o(T^{a-1})$ and $Q_{T,222}=o(T^{a-1})$ uniformly for $b>0$ as $T\rightarrow\infty$, respectively.
The proof of the lemma is complete.
\hfill$\Box$

\textbf{Lemma A2}. {\sl Let $q=q(T)=b(\frac{2}{m}\ln T)^{-1/\alpha}$
for some constant $b>0$. Under the conditions of Theorem 2.2
we have
$$P\left\{\max_{t\in \mathfrak{R}(q)\cap[0,T^{a}]}Y_{m:n}(t)> v_{T}\right\}=T^{a-1}e^{-x-r+\sqrt{2rm}\mathcal{N}}(1+o(1)),$$
as $b\downarrow0$ and $T\rightarrow\infty$;
$$P\left\{\max_{t\in \mathfrak{R}(p)\cap[0,T^{a}]}Y_{m:n}(t)> v_{T}^{*}\right\}=T^{a-1}e^{-x-r+\sqrt{2rm}\mathcal{N}}(1+o(1))$$
as $T\rightarrow\infty$, and
\begin{eqnarray*}
&&P\left\{\max_{t\in \mathfrak{R}(q)\cap[0,T^{a}]}Y_{m:n}(t)> v_{T},
\max_{t\in \mathfrak{R}(p)\cap[0,T^{a}]}Y_{m:n}(t)> v_{T}^{*}\right\}\\
&&\ \ \ \ \ =T^{a-1}\mathcal{H}_{d,m,\alpha}^{\ln
(\mathcal{H}_{m,\alpha})+x, \ln(\mathcal{H}_{d,m,\alpha})+y}e^{-r+\sqrt{2rm}\mathcal{N}}(1+o(1)),
\end{eqnarray*}
as $b\downarrow0$ and $T\rightarrow\infty$, where $v_{T}$ and $v_{T}^{*}$ are defined in (\ref{equ413}) and (\ref{equ414}), respectively. }

\textbf{Proof:} The first assertion is just the one in Lemma A1 and we present it here just for citing easily. Recall that $\mathfrak{R}(p)=\mathfrak{R}(d(\frac{2}{m}\ln
T)^{-1/\alpha})$ with $d>0$  is a Pickands grid under the conditions of Theorem 2.2.
Thus, the second assertion can be proved by repeating the proof of Proposition 2.1 and Theorem 4.1 of D\c{e}bicki, Hashorva, Ji and Tabi\'{s} (2015) by replacing the Pickands type constant $\mathcal{H}_{m,\alpha}$
by $\mathcal{H}_{d,m,\alpha}$. We omit the details.

We consider the third assertion. It can be proved quite similar to the
proof of Proposition 2.1 of D\c{e}bicki, Hashorva, Ji and Tabi\'{s} (2015), that
$$P\left(\max_{t\in [0,\lambda u^{-2/\alpha}]}Y_{m:m}(t)>u+\frac{x}{u},\max_{t\in [0,\lambda u^{-2/\alpha}]\cap \mathfrak{R}(d)}Y_{m:m}(t)>u\right)= \mathcal{H}_{d,m,\alpha}^{x,0}(\lambda)\Psi^{m}(u)(1+o(1))$$
as $u\rightarrow\infty$, where $\mathfrak{R}(d)=(dku^{-2/\alpha})$ with $k\in \mathbb{N}$ is a Pickands grid in $\mathbb{R}$ and
$$\mathcal{H}_{d,m,\alpha}^{x,0}(\lambda)=\int_{\mathbb{R}^{m}}e^{\sum_{i=1}^{m}w_{i}}
P\left(\sup_{t\in[0,\lambda]}\min_{1\leq i\leq m}(\sqrt{2}B_{\alpha/2}^{(i)}(t)-t^{\alpha}>w_{i}+x,
\sup_{kd\in[0,\lambda]}\min_{1\leq i\leq m}(\sqrt{2}B_{\alpha/2}^{(i)}(kd)-(kd)^{\alpha}>w_{i}\right)d\mathbf{w}.$$
It also can be proved in a similar way as for Proposition 2.1 and Theorem 4.1 of D\c{e}bicki, Hashorva, Ji and Tabi\'{s} (2015)  that
$$\mathcal{H}_{d,m,\alpha}^{x,0}:=\lim_{\lambda\rightarrow\infty}
\frac{\mathcal{H}_{d,m,\alpha}^{x,0}(\lambda)}{\lambda}\in
(0,\infty).$$
Now, by the same arguments as the proof of Theorem 2.1 in D\c{e}bicki, Hashorva, Ji and Ling (2014),  we  can show that
\begin{eqnarray*}
&&P\left(\max_{t\in [0,\lambda u^{-2/\alpha}]}Y_{m:n}(t)>u+\frac{x}{u},\max_{t\in [0,\lambda u^{-2/\alpha}]\cap \mathfrak{R}(d)}Y_{m:n}(t)>u\right)\\
&&=C_{n}^{m}P\left(\max_{t\in [0,\lambda u^{-2/\alpha}]}Y_{m:m}(t)>u+\frac{x}{u},\max_{t\in [0,\lambda u^{-2/\alpha}]\cap \mathfrak{R}(d)}Y_{m:m}(t)>u\right)(1+o(1))\\
&&=C_{n}^{m}\mathcal{H}_{d,m,\alpha}^{x,0}(\lambda)\Psi^{m}(u)(1+o(1)),
\end{eqnarray*}
as $u\rightarrow\infty$. Now, using the above facts and by similar arguments as for Theorem 7.1 and  Corollary 7.3  of Piterbarg
(1996), we have
\begin{eqnarray*}
\label{eq318}
&&\lim_{b\downarrow0}P\left\{\max_{t\in \mathfrak{R}(q)\cap[0,T^{a}]}Y_{m:n}(t)> v_{T}+ \frac{x}{v_{T}},
\max_{t\in \mathfrak{R}(p)\cap[0,T^{a}]}Y_{m:n}(t)>v_{T}\right\}\\
&&=P\left\{\max_{t\in [0,T^{a}]}Y_{m:n}(t)>v_{T}+ \frac{x}{v_{T}},
\max_{t\in \mathfrak{R}(p)\cap[0,T^{a}]}Y_{m:n}(t)>v_{T}\right\}\\
&&=T^{a}C_{n}^{m}\mathcal{H}_{d,m,\alpha}^{x,0}v_{T}^{2/\alpha}\Psi^{m}(v_{T})(1+o(1))
\end{eqnarray*}
as $T\rightarrow\infty$. Now, to complete the proof of the third assertion, we only need to make some transform. Using
(\ref{equ413}) and (\ref{equ414}), we get
\begin{eqnarray*}
v_{T}
&=&\frac{x+r-\sqrt{2rm}\mathcal{N}}{a_{m,T}}+b_{m,T}+o(a_{m,T}^{-1})\\
&=&v_{T}^{*}+b_{m,T}-b_{d,m,T}+(x-y)/a_{m,T}+o(a_{m,T}^{-1})\\
&=&v_{T}^{*}+\frac{\ln (\mathcal{H}_{m,\alpha})-\ln (\mathcal{H}_{d,m,\alpha})+x-y}{v_{T}^{*}}+O\left((\ln\ln (T))^{2}(\ln T)^{-3/2}\right).
\end{eqnarray*}
Observing that $v_{T}^{*}\thicksim (\frac{2}{m}\ln
T)^{1/2}$, we see that the reminder $O(\cdot)$ plays a
negligible role. Therefore, using the definition of $v_{T}^{*}$ again, we have
\begin{eqnarray*}
&&\lim_{b\downarrow0}P\left\{\max_{t\in \mathfrak{R}(q)\cap[0,T^{a}]}Y_{m:n}(t)> v_{T},
\max_{t\in \mathfrak{R}(p)\cap[0,T^{a}]}Y_{m:n}(t)>v_{T}^{*}\right\}\\
&&=T^{a}C_{n}^{m}\mathcal{H}_{d,m,\alpha}^{Z_{x,y},0}(v_{T}^{*})^{2/\alpha}\Psi(v_{T}^{*})(1+o(1))\\
&&=T^{a-1}\mathcal{H}_{d,m,\alpha}^{Z_{x,y},0}(\mathcal{H}_{d,m,\alpha})^{-1}e^{-y-r+\sqrt{2rm}\mathcal{N}}(1+o(1)),
\end{eqnarray*}
where $Z_{x,y}=\ln
(\mathcal{H}_{m,\alpha})-\ln
(\mathcal{H}_{d,m,\alpha})+x-y$.
Next, changing the variables in the definition of
$\mathcal{H}_{d,m,\alpha}^{x,y}$ we get that
$\mathcal{H}_{d,m\alpha}^{Z_{x,y},0}(\mathcal{H}_{d,m,\alpha})^{-1}e^{-y}
=\mathcal{H}_{d,m,\alpha}^{\ln
(\mathcal{H}_{m,\alpha})+x, \ln
(\mathcal{H}_{d,m,\alpha})+y}$, which completes the proof of the lemma. \hfill$\Box$

\subsection{Appendix B}

In this subsection, we give three technical lemmas which are used for the proof of Lemma 3.3. Recall that $u_{T}=
b_{m,T}+x/a_{m,T}$, $u_{T}^{*}=
b_{m,T}^{*}+y/a_{m,T}$, where
$b_{m,T}^{*}=b_{m,T}^{p}$ for  sparse grids and
$b_{m,T}^{*}=b_{d,m,T}$ for  Pickands grids, and
$r^{(h)}(t,s)=hr(t,s)+(1-h)\varrho(t,s)$ with $h\in[0,1]$.
Let
$$\theta(t,s)=\max\{|r(|t-s|)|, |\varrho(t,s)|\}$$
and
 $$\theta(z)=\sup_{0\leq s,t\leq T,|s-t|>z}\{\theta(t,s)\}.
$$
It is easy to see from assumption (\ref{eq1.1}) that for any $\varepsilon>0$,
$\theta(\varepsilon)<1$
for all sufficiently large $T$. Furthermore, choose positive constants $a,c$ be such that
\begin{eqnarray}
\label{eqTB}
0<c<a<\frac{1-\theta(\varepsilon)}{1+\theta(\varepsilon)}<1
\end{eqnarray}
 for all sufficiently large $T$ and for some $\varepsilon>0$ which will be chosen in the sequel.
It follows from (\ref{eq1.1}) again that for all
$|s-t|\leq\varepsilon<2^{-1/\alpha}$,
\begin{eqnarray}
\label{eqT20}
\frac{1}{2}|s-t|^{\alpha}\leq 1-r(s-t)\leq 2|s-t|^{\alpha},
\end{eqnarray}
which will be used in this section.\\

\textbf{Lemma B1}. {\sl Under the conditions of Lemma 3.3, we have
\begin{eqnarray}
\label{eqB11}
\sum_{t\in \mathfrak{R}(q)\cap\mathbf{O}_{i}, s\in \mathfrak{R}(q)\cap\mathbf{O}_{j},\atop t\neq s, 1\leq i,j\leq l }|r(|t-s|)-\varrho(s,t)|
\int_{0}^{1}\frac{u_{T}^{-2(m-1)}}{(1-r^{(h)}(t,s))^{m/2}}\exp\left(-\frac{mu_{T}^{2}}{1+r^{(h)}(t,s)}\right)dh
\rightarrow 0
\end{eqnarray}
as $T\rightarrow\infty$.
}

\textbf{Proof:} Recall that  $\mathfrak{R}(q)$ is a Pickands grid.
First, we consider the case that $s,t$ in the
same interval $\mathbf{O_{i}}$.
Split the sum (\ref{eqB11}) into two parts as
\begin{eqnarray}
\label{eq402}
\sum_{t\in \mathfrak{R}(q)\cap\mathbf{O}_{i}, s\in \mathfrak{R}(q)\cap\mathbf{O}_{j},\atop t\neq s, 1\leq i,j\leq l,|t-s|\leq\varepsilon  }
+\sum_{t\in \mathfrak{R}(q)\cap\mathbf{O}_{i}, s\in \mathfrak{R}(q)\cap\mathbf{O}_{j},\atop t\neq s, 1\leq i,j\leq l,|t-s|>\varepsilon  }=:J_{T,1}+J_{T,2}.
\end{eqnarray}
We deal with $J_{T,1}$ and note that in this case, we have
$\varrho(t,s)-r(|t-s|)=\rho(T)(1-r(|t-s|))$.
By  (\ref{eq1.1}) we can choose small enough $\varepsilon>0$ such that
$\varrho(t,s)=r(|t-s|)+(1-r(|t-s|))\rho(T)\sim
r(|t-s|)$ for sufficiently large $T$ and $|t-s|\leq\varepsilon$.
The definition of $u_{T}$ implies
\begin{eqnarray}
\label{eq403}
u^{2}_{T}=\frac{2}{m}\ln T+ \frac{2}{m}\ln[(\ln T)^{1/\alpha-m/2}]+O(1).
\end{eqnarray}
Consequently, we have
\begin{eqnarray*}
\label{eq404}
J_{T,1}&\leq&\mathcal{C}\sum_{t\in \mathfrak{R}(q)\cap\mathbf{O}_{i}, s\in \mathfrak{R}(q)\cap\mathbf{O}_{j},\atop t\neq s, 1\leq i,j\leq l,|t-s|\leq\varepsilon  }|r(|t-s|)-\varrho(s,t)|\frac{u_{T}^{-2(m-1)}}{(1-r(|t-s|))^{m/2}}\exp\left(-\frac{mu^{2}_{T}}{1+r(|t-s|)}\right)\nonumber\\
&= & \mathcal{C} Tb^{-1}u_{T}^{2/\alpha}\rho(T)u_{T}^{-2(m-1)}\exp\left(-\frac{mu^{2}_{T}}{2}\right)\sum_{t\in \mathfrak{R}(q)\cap[0,T],|t|\leq \varepsilon }(1-r(t))^{1-m/2}
\exp\left(-\frac{m(1-r(t))u^{2}_{T}}{2(1+r(t))}\right).
\end{eqnarray*}
Now, by (\ref{eqT20}) and (\ref{eq403}) and noting that $\rho(T)=r/\ln T=O(u_{T}^{-2})$
\begin{eqnarray*}
J_{T,1}&\leq & \mathcal{C} b^{-1}u_{T}^{-m}\sum_{ t\in \mathfrak{R}(q)\cap[0,T],|t|\leq \varepsilon }\sqrt{2}|t|^{\alpha-m\alpha/2}\exp\left(-\frac{m|t|^{\alpha}u^{2}_{T}}{8}\right)\nonumber\\
&\leq & \mathcal{C} b^{-(1-\alpha+m\alpha/2)}u_{T}^{-2}\sum_{k=1}^{\infty}(k)^{\alpha-m\alpha/2}e^{-\frac{1}{4}m(kb)^{\alpha}}\nonumber\\
&\leq & \mathcal{C} b^{-(1-\alpha+m\alpha/2)}u_{T}^{-2},
\end{eqnarray*}
which shows $J_{T,1}\rightarrow0$ uniformly for $b>0$ as $T\rightarrow\infty$.\\
Using the fact that $u_{T}\thicksim (\frac{2}{m}\ln
T)^{1/2}$, we obtain
\begin{eqnarray}
\label{eq405}
J_{T,2}&\leq&\mathcal{C}u_{T}^{-2(m-1)}\sum_{t\in \mathfrak{R}(q)\cap\mathbf{O}_{i}, s\in \mathfrak{R}(q)\cap\mathbf{O}_{j},\atop t\neq s, 1\leq i,j\leq l,|t-s|>\varepsilon  }
\exp\left(-\frac{mu^{2}_{T}}{1+r(|t-s|)}\right)\nonumber\\
&\leq & \mathcal{C}T^{1+a}b^{-2}u_{T}^{4/\alpha}u_{T}^{-2(m-1)}
     \exp\left(-\frac{mu^{2}_{T}}{1+\theta(\varepsilon)}\right)\nonumber\\
&\leq & \mathcal{C} T^{1+a}b^{-2}u_{T}^{4/\alpha}u_{T}^{-2(m-1)}(T)^{-\frac{2}{1+\theta(\varepsilon)}}\nonumber\\
&\leq & \mathcal{C} T^{a-\frac{1-\theta(\varepsilon)}{1+\theta(\varepsilon)}}b^{-2}(\ln T)^{2/\alpha-m+1}.
\end{eqnarray}
Thus, $J_{T,2}\rightarrow0$ uniformly for $b>0$ as $T\rightarrow\infty$, since $a<\frac{1-\theta(\varepsilon)}{1+\theta(\varepsilon)}$.

Second, we consider the case that $t\in\mathbf{O}_{i}$ and $s\in \mathbf{O}_{j}$ for $i\neq j$.
Note that in this case, $|t-s|\geq T^{c}$ and $\varrho(s,t)=\rho(T)$ for  $s\in \mathbf{O}_{i}$
and $t\in \mathbf{O}_{j}$, $i\neq j$. Choose $\beta$ such that
$0<c<a<\beta<\frac{1-\theta(\varepsilon)}{1+\theta(\varepsilon)}$ and
split the sum (\ref{eqB11}) into two parts as
\begin{eqnarray}
\label{eq402}
\sum_{t\in \mathfrak{R}(q)\cap\mathbf{O}_{i}, s\in \mathfrak{R}(q)\cap\mathbf{O}_{j},\atop  1\leq i\neq j\leq l,|t-s|\leq T^{\beta}  }
+\sum_{t\in \mathfrak{R}(q)\cap\mathbf{O}_{i}, s\in \mathfrak{R}(q)\cap\mathbf{O}_{j},\atop  1\leq i\neq j\leq l,|t-s|> T^{\beta}  }=:S_{T,1}+S_{T,2}.
\end{eqnarray}
For $S_{T,1}$, with the similar derivation as for
(\ref{eq405}), we have
\begin{eqnarray}
\label{eq407}
S_{T,1}&\leq&\mathcal{C}u_{T}^{-2(m-1)}\sum_{t\in \mathfrak{R}(q)\cap\mathbf{O}_{i}, s\in \mathfrak{R}(q)\cap\mathbf{O}_{j},\atop  1\leq i\neq j\leq l,|t-s|\leq T^{\beta}  }
\exp\left(-\frac{mu^{2}_{T}}{1+r(|t-s|)}\right)\nonumber\\
&\leq & \mathcal{C}T^{1+\beta}b^{-2}u_{T}^{4/\alpha}u_{T}^{-2(m-1)}
     \exp\left(-\frac{mu^{2}_{T}}{1+\theta(\varepsilon)}\right)\nonumber\\
&\leq & \mathcal{C} T^{1+\beta}b^{-2}u_{T}^{4/\alpha}u_{T}^{-2(m-1)}(T)^{-\frac{2}{1+\theta(\varepsilon)}}\nonumber\\
&\leq & C T^{\beta-\frac{1-\theta(\varepsilon)}{1+\theta(\varepsilon)}}b^{-2}(\ln T)^{2/\alpha-m+1}.
\end{eqnarray}
Consequently, since $\beta<\frac{1-\theta(\varepsilon)}{1+\theta(\varepsilon)}$, we have $S_{T,1}\rightarrow0$ uniformly for $b>0$ as $T\rightarrow\infty$.\\
For $S_{T,2}$, we need more precise estimation.
By condition (\ref{eq1.2}), there exist constants $C>0$ and $K>0$
such that
$$\theta(t)\ln t\leqslant K$$
for $T$ sufficiently large and $t$ satisfying
$t\geq C$. Thus for $T$ large enough and for $t$ such
that $t\geq T^\beta$,
$\theta(t)\leq K/\ln(T^{\beta})$.
Now making use of (\ref{eq403}), we
obtain
\begin{eqnarray}
\label{Tan2}
&&T^{2}u_{T}^{4/\alpha}(\ln T)^{-m}\exp\left(-\frac{mu_{T}^{2}}{1+\theta(T^{\beta})}\right)\nonumber\\
&&\leq T^{2}u_{T}^{4/\alpha}(\ln T)^{-m}\exp\left(-\frac{mu_{T}^{2}}{1+K/\ln(T^{\beta})}\right)\nonumber\\
&&\leq O(1) T^{2}(\ln T)^{2/\alpha-m}\left(T^{-2}(\ln T)^{-(2/\alpha-m)}\right)^{\frac{1}{1+K/\ln(T^{\beta})}}\nonumber\\
&&\leq O(1)T^{(2K/\ln(T^{\beta}))/(1+K/\ln(T^{\beta}))}(\ln T)^{((2/\alpha-m)K/\ln (T)^{\beta})/(1+K/\ln (T^{\beta}))}\nonumber\\
&&=O(1).
\end{eqnarray}
Therefore, by a similar argument as for the proof of Lemma 6.4.1 of
Leadbetter et al. (1983) we obtain
\begin{eqnarray}
\label{eq409}
S_{T,2}&\leq&\mathcal{C}u_{T}^{-2(m-1)}\sum_{t\in \mathfrak{R}(q)\cap\mathbf{O}_{i}, s\in \mathfrak{R}(q)\cap\mathbf{O}_{j},\atop  1\leq i\neq j\leq l,|t-s|> T^{\beta}  }|r(|t-s|)-\rho(T)|\exp\left(-\frac{mu^{2}_{T}}{1+\theta(T^{\beta})}\right)\nonumber\\
&\leq&\mathcal{C}Tb^{-1}u_{T}^{2/\alpha}u_{T}^{-2(m-1)}\exp\left(-\frac{mu^{2}_{T}}{1+\theta(T^{\beta})}\right)
      \sum_{t\in \mathfrak{R}(q)\cap[0,T], t>T^{\beta}}|r(t)-\rho(T)|\nonumber\\
&=& \mathcal{C}T^{2}(\ln T)^{-m}u_{T}^{4/\alpha}\exp\left(-\frac{mu^{2}_{T}}{1+\theta(T^{\beta})}\right)\cdot b^{-1}\frac{\ln T}{Tu_{T}^{2/\alpha}}
     \sum_{t\in \mathfrak{R}(q)\cap[0,T], t>T^{\beta}}|r(t)-\rho(T)|\nonumber\\
&\leq& \mathcal{C}b^{-1}\frac{\ln T}{Tu_{T}^{2/\alpha}}
     \sum_{t\in \mathfrak{R}(q)\cap[0,T], t>T^{\beta}}|r(t)-\rho(T)|\nonumber\\
&\leq& \mathcal{C}b^{-1}\frac{1}{\beta Tu_{T}^{2/\alpha}}
     \sum_{t\in \mathfrak{R}(q)\cap[0,T], t>T^{\beta}}|r(t)\ln t-r|+ \mathcal{C}b^{-1}\frac{r}{Tu_{T}^{2/\alpha}}
     \sum_{t\in \mathfrak{R}(q)\cap[0,T], t>T^{\beta}}\big|1-\frac{\ln T}{\ln t}\big|.
\end{eqnarray}
By condition (\ref{eq1.2}) the first term on the
right-hand-side of (\ref{eq409}) tends to 0 uniformly for $b>0$ as $T\to
\infty$. Furthermore, the second term of the right-hand-side of
(\ref{eq409}) also tends to 0 by an integral estimate as follows (
see also the proof of Lemma 6.4.1 of Leadbetter et al.\ (1983))
\begin{eqnarray*}
\mathcal{C}b^{-1}\frac{r}{Tu_{T}^{2/\alpha}}
     \sum_{t\in \mathfrak{R}(q)\cap[0,T], t>T^{\beta}}\big|1-\frac{\ln T}{\ln t}\big|
&\leq& \mathcal{C}b^{-1}\frac{r}{Tu_{T}^{2/\alpha}}\frac{1}{\ln T^{\beta}}
     \sum_{t\in \mathfrak{R}(q)\cap[0,T], t>T^{\beta}}\big|\ln t-\ln T\big|\\
&=& O\left(b^{-1}\frac{r}{\ln T^{\beta}}\int_0^1 |\ln x|dx\right),
\end{eqnarray*}
which shows that $S_{T,2}\rightarrow 0$ uniformly for $b>0$ as $T\rightarrow\infty$.
The proof of the lemma is complete.\hfill$\Box$\\

\textbf{Lemma B2}. {\sl Under the conditions of Lemma 3.3, we have
\begin{eqnarray}
\label{B2}
\sum_{t\in \mathfrak{R}(p)\cap\mathbf{O}_{i}, s\in \mathfrak{R}(p)\cap\mathbf{O}_{j},\atop t\neq s, 1\leq i,j\leq l }|r(|t-s|)-\varrho(s,t)|
\int_{0}^{1}\frac{u_{T}^{-2(m-1)}}{(1-r^{(h)}(t,s))^{m/2}}\exp\left(-\frac{m(u_{T}^{*})^{2}}{1+r^{(h)}(t,s)}\right)dh
\rightarrow 0
\end{eqnarray}
as $T\rightarrow\infty$.
}

\textbf{Proof:}
The proof is the same as that of Lemma B1, we omit the details.\\

\textbf{Lemma B3}. {\sl Under the conditions of Lemma 3.3, we have
\begin{eqnarray}
\label{eqB30}
\sum_{t\in \mathfrak{R}(q)\cap\mathbf{O}_{i}, s\in \mathfrak{R}(p)\cap\mathbf{O}_{j},\atop t\neq s, 1\leq i,j\leq l }|r(|t-s|)-\varrho(s,t)|
\int_{0}^{1}\frac{u_{T}^{-2(m-1)}}{(1-r^{(h)}(t,s))^{m/2}}\exp\left(-\frac{\frac{1}{2}m[u_{T}^{2}+(u_{T}^{*})^{2}]}{1+r^{(h)}(t,s)}\right)dh
\rightarrow 0
\end{eqnarray}
as $T\rightarrow\infty$.
}

\textbf{Proof:} Recall that $\mathfrak{R}(p)$ can be a sparse grid or Pickands grid. We only show the case that $\mathfrak{R}(p)$ is a sparse grid,
since the proof of the remaining case is similar.
First, we consider the case that $s,t$ in the
same interval $\mathbf{O_{i}}$.
Split the sum (\ref{eqB30}) into two parts as
\begin{eqnarray}
\label{eqB31}
\sum_{t\in \mathfrak{R}(q)\cap\mathbf{O}_{i}, s\in \mathfrak{R}(p)\cap\mathbf{O}_{j},\atop t\neq s, 1\leq i=j\leq l, |t-s|\leq\varepsilon }
+\sum_{t\in \mathfrak{R}(q)\cap\mathbf{O}_{i}, s\in \mathfrak{R}(p)\cap\mathbf{O}_{j},\atop t\neq s, 1\leq i=j\leq l, |t-s|>\varepsilon }=:W_{T,1}+W_{T,2}.
\end{eqnarray}
We deal with $W_{\mathbf{T},1}$ and note that in this case, by the definition of the field $\xi_{T}(t)$, we have
$\varrho(t,s)-r(|t-s|)=\rho(T)(1-r(|t-s|))$.
By  (\ref{eq1.1}) we can choose small enough $\varepsilon>0$ such that
$\varrho(t,s)=r(|t-s|)+(1-r(|t-s|))\rho(T)\sim
r(|t-s|)$ for sufficiently large $T$ and $|t-s|\leq\varepsilon$.
By the definitions of $u_{T}$ and $u_{T}^{*}$, we have
\begin{eqnarray}
\label{eqB32}
w_{T}^{2}:=\frac{1}{2}(u_{T}^{2}+(u_{T}^{*})^{2})=\frac{2}{m}\ln T+\frac{1}{m}\ln (\ln T)^{1/\alpha-m/2}+\frac{1}{m}\ln (p^{-1}(\ln T)^{-m/2})+O(1).
\end{eqnarray}
Consequently, we have
\begin{eqnarray*}
\label{eq404}
W_{T,1}&\leq&\mathcal{C}u_{T}^{-2(m-1)}\sum_{t\in \mathfrak{R}(q)\cap\mathbf{O}_{i}, s\in \mathfrak{R}(p)\cap\mathbf{O}_{j},\atop t\neq s, 1\leq i=j\leq l, |t-s|\leq\varepsilon }
|r(|t-s|)-\varrho(s,t)|\frac{1}{(1-r(|t-s|))^{m/2}}\exp\left(-\frac{mw^{2}_{T}}{1+r(|t-s|)}\right)\nonumber\\
&\leq &\mathcal{ C}Tp^{-1}u_{T}^{-2(m-1)}\rho(T)\exp\left(-\frac{mw^{2}_{T}}{2}\right)\sum_{t\in \mathfrak{R}(q)\cap[0,T], |t|\leq \varepsilon }(1-r(t))^{1-m/2}
\exp\left(-\frac{m(1-r(t))w^{2}_{T}}{2(1+r(t))}\right),
\end{eqnarray*}
then by (\ref{eqT20}) and (\ref{eqB32})
\begin{eqnarray*}
W_{T,1}
&\leq & \mathcal{C}p^{-1/2}u_{T}^{-2m}(\ln T)^{-1/2\alpha+m/2} \sum_{ t\in \mathfrak{R}(q)\cap[0,T],|t|\leq \varepsilon }|t|^{\alpha-m\alpha/2}\exp\left(-\frac{m|t|^{\alpha}w^{2}_{T}}{8}\right)\nonumber\\
&\leq & \mathcal{C}(b)^{\alpha-m\alpha/2} [p(\ln T)^{1/\alpha}]^{-1/2}u_{T}^{-2}\sum_{k=1}^{\infty}(k)^{\alpha-m\alpha/2}e^{-\frac{1}{4}(kb)^{\alpha}}\nonumber\\
&\leq & \mathcal{C} [p(\ln T)^{1/\alpha}]^{-1/2}u_{T}^{-2},
\end{eqnarray*}
which shows $W_{T,1}\rightarrow0$ uniformly for $b>0$ as $T\rightarrow\infty$, since $p(\ln T)^{1/\alpha}\rightarrow\infty$ for sparse grid.\\
Using the fact that $w_{T}\thicksim (\frac{2}{m}\ln
T)^{1/2}$, we obtain
\begin{eqnarray*}
\label{eq405}
W_{T,2}&\leq&\mathcal{C}\sum_{t\in \mathfrak{R}(q)\cap\mathbf{O}_{i}, s\in \mathfrak{R}(p)\cap\mathbf{O}_{j},\atop t\neq s, 1\leq i=j\leq l, |t-s|\leq\varepsilon }
\exp\left(-\frac{w^{2}_{T}}{1+r(|t-s|)}\right)\nonumber\\
&\leq & \mathcal{C}T^{1+a}u_{T}^{2/\alpha}p^{-1}b^{-1}u_{T}^{-2(m-1)}
     \exp\left(-\frac{w^{2}_{T}}{1+\theta(\varepsilon)}\right)\nonumber\\
&\leq & \mathcal{C} T^{1+a}u_{T}^{2/\alpha}p^{-1}b^{-1}u_{T}^{-2(m-1)}(T)^{-\frac{2}{1+\theta(\varepsilon)}}\nonumber\\
&\leq & C T^{a-\frac{1-\theta(\varepsilon)}{1+\theta(\varepsilon)}}b^{-1}(\ln T)^{2/\alpha-m+1}.
\end{eqnarray*}
Thus, $W_{T,2}\rightarrow0$ uniformly for $b>0$ as $T\rightarrow\infty$, since $a<\frac{1-\theta(\varepsilon)}{1+\theta(\varepsilon)}$.

Second, we consider the case that $t,s$ in the
different intervals $\mathbf{O}_{i}$ and $\mathbf{O}_{j}$ for $i\neq j$. Since the proof is similar as that of Lemma B1, we omit the details.
\hfill$\Box$

\subsection{Appendix C}
In this subsection, we give a comparison inequality, which is a slight extension of Theorem 2.4 of D\c{e}bicki, Hashorva, Ji and Ling (2017).

Denote by
$\mathcal{X}=(X_{il})_{d\times n}$ and $\mathcal{Y}=(Y_{il})_{d\times n}$  two random arrays with $N(0,1)$ components, and let
$\Sigma^{(1)}=(\sigma^{(1)}_{il,jk})_{{dn \times dn}}$ and $\Sigma^{(0)}=(\sigma^{(0)}_{il,jk})_{dn \times dn}$ be the covariance matrices of $\mathcal{X}$ and $\mathcal{Y}$, respectively, with
$ \sigma^{(1)}_{il,jk} := E{X_{il}X_{jk}}$ and $\sigma^{(0)}_{il,jk} := E{Y_{il}Y_{jk}}, 1\le i,j\le d, 1\le l,k\le n.$
Furthermore, define $\mathbf{X}_{(m)}=(X_{1(m)}, \cdots, X_{d(m)}), 1\le r\le n$ to be the  $m$-th order statistics vector generated by $\mathcal{X}$ as follows
\begin{eqnarray*}
X_{i(1)}= \min_{1\le l\le n} X_{il} \le  \cdots \le X_{i(m)}\le\cdots \le \max_{1\le l\le n}X_{il}= X_{i(n)},\quad 1\le i\le d.
\end{eqnarray*}
Similarly, we write    $\mathbf{Y}_{(m)}
=(Y_{1(m)}, \cdots, Y_{d(m)})$ which is generated by $\mathcal{Y}$.
Assume that  the columns of both $\mathcal{X}$ and $\mathcal{Y}$ are mutually independent, i.e.,
\begin{eqnarray*}
\sigma_{il,jk}^{(\kappa)}= \sigma_{ij}^{(\kappa)}\mathbb I\{l=k\}, \quad  1\le i,  j\le d, 1\le l, k\le n, \kappa=0,1,
\end{eqnarray*}
with some $\sigma_{ij}^{(\kappa)},  1\le i,  j\le d, \kappa=0,1$, where $\mathbb I\{\cdot\}$ stands for the indicator function.

\textbf{Lemma C1}. {\sl Let $\mathbf{u_{T}}=(u_{T1},\ldots,u_{Td})$. Suppose that  $u_{Ti}\rightarrow \infty $
and $u_{Ti}/u_{Tj}\rightarrow 1$ for all $i,j\in \{1,2,\ldots,d\}$ as $T_{i}\rightarrow\infty$. Then for sufficiently large $\mathbf{T}$ , we have
for some constant $\mathcal{C}>0$
\begin{eqnarray*}
&&|P(\mathbf{X}_{(m)}\leq \mathbf{u_{T}})-P(\mathbf{Y}_{(m)}\leq \mathbf{u_{T}})|\\
&&\leq \mathcal{C} u_{T1}^{-2(n-m)}\sum_{1\leq i<j\leq d}|\sigma_{ij}^{(0)}-\sigma_{ij}^{(1)}|\int_{0}^{1}(1-\delta_{ij}^{(h)})^{-(n-m+1)/2}
\exp\left(-\frac{(n-m+1)(u_{Ti}^{2}+u_{Tj}^{2})}{2(1+|\delta_{ij}^{(h)}|)}\right)dh,
\end{eqnarray*}
where $\delta_{ij}^{(h)}=h\sigma_{ij}^{(0)}+(1-h)\sigma_{ij}^{(1)}$.
}

\textbf{Proof:} The proof is similar to that of Theorem 2.4 of  D\c{e}bicki, Hashorva, Ji and Ling (2017) with some changes.
Let $(Z_i, Z_j)$ be a bivariate standard normal random vector with correlation $|\delta_{ij}^{(h)}|$.  By a
similar argument as the proof in P225 of Leadbetter et al. (1983), we can show (for large $T_{i}$)
\begin{eqnarray*}
&&P\{Z_i >u_{Ti}, Z_j >u_{Tj}\} \leq \frac{\mathcal{C}}{u_{Ti}(u_{Tj}-|\delta_{ij}^{(h)}|u_{Ti})}\phi(u_{Ti}, u_{Tj}; |\delta_{il}^{(h)}|)
\leq \frac{\mathcal{C}}{u_{T1}^{2}}\phi(u_{Ti}, u_{Tj}; |\delta_{ij}^{(h)}|),
\end{eqnarray*}
where $\phi(u,v,r)$ is the probability density function of a two dimensional normal random variables.  Using the above
inequality to replace (4.28) in the proof of Theorem 2.4 of D\c{e}bicki, Hashorva, Ji and Ling (2017), we can prove Lemma C1.\hfill$\Box$

\COM{\bigskip
{\bf Acknowledgement}: We would like to thank ...}


\begin{thebibliography}{100} \small

\bibitem{} Albin JPM. (1990) On extremal theory for stationary processes. Ann Probab., 18:92-128











\bibitem{} D\c{e}bicki, K., Hashorva, E., Ji, L. and Ling, C. (2015) Extremes of order statistics of stationary processes. Test, 24, 229-248.

\bibitem{} D\c{e}bicki, K., Hashorva, E., Ji, L. and Ling, C. (2017) Comparision inequality for order statistics of Gaussian arrays.
Latin American Journal of Probability and Mathematical Statistics, 14, 93-116.

\bibitem{} D\c{e}bicki, K., Hashorva, E., Ji, L. and Tabi\'{s}, K. (2014) On the probability of conjunctions of stationary Gaussian processes. Statist.
Probab. Lett., 88, 141-148.

\bibitem{} D\c{e}bicki, K., Hashorva, E., Ji, L. and Tabi\'{s}, K. (2015) Extremes of vector-valued Gaussian processes: Exact asymptotics. Stochastic Process. Appl., 125(11), 4039-4065.


\bibitem{} D\c{e}bicki, K., Kosinski, K.M. (2018) An Erd\"{o}s-R\'{e}v\'{e}sz type law of the iterated logarithm for order statistics of a stationary Gaussian process. J. Theor. Probab.,  31(1), 579-597.


\bibitem{} Hashorva, E., Tan, Z. (2015)  Piterbarg's max-discretisation theorem for stationary
vector Gaussian processes observed on different grids. Statistics, 49(2), 338-360.



\bibitem{} H\"{u}sler J. (2004) Dependence between extreme values of discrete and continuous time locally stationary Gaussian processes. Extremes, 7, 179-190.

\bibitem{} H\"{u}sler, J., Piterbarg, V.I. (2004) Limit theorem for maximum of the storage process with
fractional Brownian motion as input. Stochastic Process. Appl., 114,
231-250.

\bibitem{} Leadbetter, M.R., Lindgren, G. (1982) Extreme value theory
for continuous parameter stationary processes. Z. Wahrsch. verw. Gebiete. 60, 1-20.


\bibitem{} Leadbetter, M.R., Lindgren, G. and Rootz\'{e}n, H., Extremes and Related Properties of Random Sequences and
Processes. Series in Statistics, Springer, New York, 1983.

\bibitem{} Ling, C., Tan, Z. (2016) On maxima of chi-processes over threshold dependent grids, Statistics, 50, 579-595.

\bibitem{} Ling, C., Peng, Z., Tan, Z. (2017) Extremes on Continuous-discretization time of Stationary Processes, Submitted for publication.

\bibitem{}  Mittal, Y., Ylvisaker, D. (1975) Limit distributions for the maxima of stationary Gaussian processes. Stochastic Processes Appl.,
3:1-18.


\bibitem{} Piterbarg, V.I. (2004) Discrete and continuous time extremes of Gaussian processes, Extremes,  7, 161-177.

\bibitem{} Piterbarg, V.I. Asymptotic Methods in the Theory of Gaussian Processes and Fields, AMS, Providence, 1996.






\bibitem{} Tan, Z., Hashorva, E., (2014) On Piterbarg's max-discretisation
theorem for multivariate stationary Gaussian processes, Journal of Mathematical Analysis and Application, 409, 299-314.

\bibitem{} Tan, Z., Wang, K., (2015) On Piterbarg's max-discretisation theorem for homogeneous Gaussian random fields,
Journal of Mathematical Analysis and Application, 429, 969-994.






\bibitem{} Turkman, K.F., Discrete and continuous time series extremes of stationary processes, in
Handbook of Statistics, T.S. Rao, S.S. Rao, and C.R. Rao, eds., Vol.
30, Time Series Methods and Aplications, Elsevier, North Holland,
2012, pp. 565-580.

\bibitem{} Zhao, C. (2018) Extremes of order statistics of stationary Gaussian processes,
Probability and  Mathematical statistics. 38(1), 61-75

\end{thebibliography}
\end{document}